\documentclass[12pt,twoside]{article}
\usepackage{amsmath,amssymb,amscd,amsthm,mathpple,mathrsfs,lcs}
\usepackage[a4paper]{geometry}
\usepackage[all]{xy}
\usepackage[english]{babel}
\usepackage[hypertex]{hyperref}
\title{Lifting systems of Galois representations
associated to abelian varieties}
\author{Rutger Noot%
\thanks{Partially supported by the European Marie Curie Research
Training Network contract MRTN-CT2003-504917}}
\date{2nd October 2006}
\begin{document}
\selectlanguage{english}
\maketitle
\begingroup
\renewcommand{\baselinestretch}{1}
\begin{abstract}
This paper treats what we call `weak geometric liftings' of Galois
representations associated to abelian varieties. 
This notion can be seen as a generalization of the idea of lifting
a Galois representation along an isogeny of algebraic groups. 
The weaker notion only takes into account an isogeny of the derived
groups and disregards the centres of the groups in question. 
The weakly lifted representations are required to be geometric in the
sense of a conjecture of Fontaine and Mazur. 
The conjecture in question states that any irreducible geometric
representation is a twist of a subquotient of an \'etale cohomology
group of an algebraic variety over a number field.

It is shown that a Galois representation associated to an abelian
variety admits a weak geometric lift to a group with simply connected
derived group. 
In certain cases, such a weak geometric lift is itself associated to
an abelian variety.
This means that the conjecture of Fontaine and Mazur is confirmed for
these representations.
In other cases, one may find a lift which can not be found back in
the \'etale cohomology of any abelian variety. 
The Fontaine--Mazur conjecture remains open for these
representations. 
Nevertheless, certain consequences of the conjecture can be
established. 
\medskip\par\noindent
\emph{2000 Mathematics Subject Classification} 11G10, 11F80, 14K15
\end{abstract}
\endgroup

\section*{Introduction}
This paper is motivated by a conjecture of Fontaine and Mazur, 
conjecture~1 of \cite{fonmaz} and  generalizes the observations
made in \cite{noot:lifting}. 
The conjecture in question aims to characterize the Tate twists of the
irreducible subquotients of the Galois representations arising from
the \'etale cohomology of algebraic varieties over number fields. 

To explain the conjecture, fix a number field $F$ and let $\cG_F$ be
the absolute Galois group of $F$. 
In \cite[\S1]{fonmaz}, 
a representation of $\cG_F$ on a finite dimensional $\Qp$-vector space
$V_p$ is called \emph{geometric} if
it is unramified outside a finite set of non-archimedean places of $F$
and if,
for each non-archimedean valuation $\bar v$ of $\bar F$ of residue
characteristic $p$, the restriction to the inertia
group at $\bar v$ is potentially semi-stable in the sense of
Fontaine, see also~\ref{intro_galois-lifts}.
An irreducible representation of $\cG_F$ on $V_p$ is said to
\emph{come from 
algebraic geometry} if it is isomorphic to a subquotient of the Galois
representation on a cohomology group $\rH_\et^i(X_{\bar F},\Qp)(m)$
for some smooth and proper $F$-scheme $X$ and integers $i$ and $m$. 
The Tate twist $(m)$ has the effect of multiplying the action of
$\cG_F$ on $\rH_\et^i(X_{\bar F},\Qp)$ by the $m$th power of the
cyclotomic character.
Fontaine and Mazur have expressed the conjecture that an irreducible
$\Qp$-linear representation $V_p$ of $\cG_F$ is geometric if and only
if it comes from algebraic geometry. 

The `if'-part of the conjecture being resolved (see \cite{tsuji:cst}),
this paper is concerned with some reflexions on the implication `only
if' in the conjecture. 
We will actually investigate a particular type of geometric
representations. 

The construction of these representations starts with the choice of an
abelian variety $A$ over $F$. 
For such an abelian variety, the representation of
$\cG_F$ on $\rH_\et^1(A_{\bar F},\Qp)$ factors through a map
\[\rho_{A,p}\colon\cG_F\rightarrow G_A(\Qp),\] where $G_A$ is the
Mum\-ford--Tate group of $A$.
After reading \cite{wintenberger:relevement}, one is tempted
to look for an isogeny $\pi\colon\wtG\rightarrow G_A$ and a map
$\trho_p\colon\cG_F\rightarrow\wtG(\Qp)$ such that
$\rho_{A,p}=\pi\circ\trho_p$ and such that $\trho_p$ defines a
geometric representation of $\cG_F$ on $\wtV_p$ for any $\Qp$-linear
representation $\wtV_p$ of $\wtG_{/\Qp}$. 

This is not exactly the point of view adopted in this paper. 
As the conjecture of Fontaine--Mazur for representations with
abelian image is quite well understood, cf.~\cite[\S6]{fonmaz}, 
the centres of the groups occurring above can be considered to
be less interesting than the derived groups. 
For this reason, given the representation
$\rho_{A,p}\colon\cG_F\rightarrow G_A(\Qp)$, we will search for a
linear algebraic group 
$\wtG$, an isogeny $\pi^\der\colon\wtG^\der\rightarrow G_A^\der$
and a map $\trho_p\colon\cG_F\rightarrow\wtG(\Qp)$ such that 
$\pi^\ad\circ\rho_{A,p}=\tilde\pi^\ad\circ\trho_p$. 
Here $\pi^\ad\colon G_A\rightarrow G_A^\ad$ and 
$\tilde\pi^\ad\colon\wtG\rightarrow\wtG^\ad$ are the natural projections and
the groups $G_A^\ad$ and $\wtG^\ad$ have been identified using the map
induced by $\pi^\der$.
If such a representation $\trho_p$ defines a geometric representation
of $\cG_F$ on any 
$\Qp$-linear representation $\wtV_p$ of $\wtG_{/\Qp}$, then we will
say that it is a \emph{weak geometric lift} of $\rho_{A,p}$. 
An advantage of working with weak geometric liftings is that,
by `correcting' liftings using characters, it is easier to produce
geometric representations than if one considers liftings along isogenies.

It turns out that for any abelian variety $A/F$, there exist a finite
extension $F'/F$, a group $\wtG$ such that $\wtG^\der$ is the
universal cover of $G_A$ and, for every prime number $p$, a weak
geometric lift $\trho_p\colon\cG_{F'}\rightarrow\wtG(\Qp)$ of the
restriction of $\rho_{A,p}$ to $\cG_{F'}$.
We refer to corollary~\ref{simply_connected_wgls_exist} for a precise
statement. 
Moreover, it follows from corollary~\ref{LADHs_dominate_all} that
every weak geometric 
lift of a $\rho_{A,p}$ is dominated by a weak geometric lift to a
group $\wtG$ with $\wtG^\der$ simply connected.
We can even conveniently `normalize' the group $\wtG$. 

It is natural to ask if the conjecture of Fontaine and Mazur is true
for any weak geometric lift of a representation associated to an
abelian variety. 
This question is not answered in this paper, but a number of
partial results are obtained. 

In section~\ref{motivic_galois_liftings} we prove the following
results. 
By combining proposition~\ref{weak_MT-lift_motive} and
remark~\ref{B'_in_thm_exists}, it follows that 
for every abelian variety $A/F$,  there
exists a group $\wtG$, a map $\wtG^\der\rightarrow G_A^\der$ and,
after replacing $F$ by a finite extension, a system $(\trho_p)$ of
weak geometric liftings of the $\rho_{A,p}$ such that 
\begin{itemize}
\item
the group $\wtG^\der$ is `not far' from the universal cover of
$G_A^\der$ and 
\item
for any representation $\wtV$ of $\wtG$, the system of
representations of $\cG_F$ on the $\wtV\otimes\Qp$ is isomorphic to
the system of $p$-adic representation associated to an abelian
variety. 
\end{itemize}
To be more precise where the first property is concerned, it means
in the first place that the group $\wtG^\der_{/\CC}$ is the product of
its simple factors $\wtG_i$. 
Secondly, it follows from well-known facts on the Mumford--Tate groups
of abelian varieties that these factors are all of classical type
($A$, $B$, $C$ or $D$). 
For the $\wtG_i$ which are of type $A_k$, $B_k$ or $C_k$, 
being `not far' from the universal cover means that they are simply 
connected. 
Where the factors $\wtG_i$ of type $D_k$ are concerned, the condition
is more difficult to state.
We have to distinguish
two subtypes, $D_k^\RR$ and $D_k^\HH$, see~\ref{D_k-factors} for the
definitions. 
The $\wtG_i$ which are of type $D_k^\RR$ in this classification are
also simply connected. 
The factors $\wtG_i$ of type $D_k^\HH$ are \emph{$h$-maximal} in the sense
of~\ref{D_k-factors}, which means that every such $\wtG_i$ is a
quotient of its universal cover by a subgroup of order $2$. 
Over $\RR$, these $h$-maximal groups are orthogonal groups.

By theorem~\ref{thm_non_spinorial_lift_from_AM}, 
the above system $(\trho_p)$ is maximal in the following sense. 
For any group $G'$, any isogeny $G^{\prime\der}\rightarrow G_A^\der$
such that $G^{\prime\der}$ is a quotient of $\wtG^\der$ and any weak
geometric lift $\rho'_p\colon\cG_F\rightarrow G'(\Qp)$ of $\rho_{A,p}$
(allowing a finite extension of $F$), the representation
$\rho'_p$ belongs to the tannakian category generated by the $p$-adic
Galois representations associated to abelian varieties and the
representations with finite image. 

The above statements about Galois representations follow from analog
properties of Hodge structures associated to abelian varieties. 
If $A_{/\CC}$ is an abelian variety over $\CC$, then the Hodge
structure on the first Betti cohomology group $\rH_\B^1(A_{/\CC}(\CC),\QQ)$
is determined by a morphism $h\colon S\rightarrow G_{A/\RR}$, where
$S=\CC^\times$ as algebraic groups over $\RR$. 
It is show in section~\ref{sect_MT_liftings} that the group $\wtG$
above is actually the Mumford--Tate group of an abelian variety
$B_{/\CC}$ and that the Hodge structure of $B$ corresponds to a
morphism $\tih\colon S\rightarrow \wtG_{/\RR}$ such that $h$ and
$\tih$ have the same projection to $G_{A/\RR}^\ad=\wtG_{/\RR}^\ad$. 
To underscore the analogy with the construction of weak geometric
liftings of Galois representations, 
such an abelian variety $B$ will be called a \emph{weak Mumford--Tate
lift} of $A$, cf.~\ref{def_MT-lift}.
This notion is very close to the notion of $A_{/\CC}$ and $B_{/\CC}$ being
`adjoint-isogenous' in the sense of definition~6.1 of Vasiu's e-print
\cite{vasiu:M-T_conj} and our variety $B_{/\CC}$ should be equal to
the variety obtained by Vasiu's `shifting process', loc.\ cit.~6.4. 

In section~\ref{sect_descent_MT-adjoint} we prove some properties
comparing the fields of definition of an abelian variety and a weak
Mumford--Tate lift. 
Using the theory of absolute Hodge motives, the above statements
concerning Galois representations are then derived from the
corresponding statements about the Hodge structures associated to weak
Mumford--Tate liftings. 
The arguments are based on those used in \cite{noot:lifting} and
\cite{paugam:good_reduction}.

In the above statements, the abelian varieties with Mumford--Tate
group of type $D_k^\HH$ stand out as an exception to the general
situation. 
These varieties deserve a separate treatment, and this is
the subject of section~\ref{AV-D_k^H-section}.
Let $F$ be a number field as above, $A$ an abelian variety over $F$ and
assume that, for the Mumford--Tate group $G_A$, the derived group
$G_{A/\CC}^\der$ is isomorphic to a product $\prod G_{A,i}$ of
$h$-maximal groups of type $D_k^\HH$. 
It can be shown that such abelian varieties do indeed exist. 
It turns out that the associated system $(\rho_{A,p})$ of $p$-adic
Galois representations admits a system of non-trivial weak geometric
liftings. 
The corollary~\ref{constr_LADH} states that there exists a group
$\wtG$ such that $\wtG^\der\rightarrow G_A^\der$ is the universal
cover and, as always after replacing $F$ by a finite extension, a
system $(\trho_p)$ of weak geometric liftings of the $\rho_{A,p}$. 
These representations are called the $p$-adic Galois representations
of \emph{lifted abelian $D_k^\HH$-type}. 
The construction makes use, again, of a lifting 
$\tih\colon S\rightarrow\wtG_{/\RR}$ of 
$h\colon S\rightarrow G_{A/\RR}$.
This time, the existence of
the system $(\trho_p)$ follows from \cite{wintenberger:relevement}.

When considering weak geometric liftings of Galois representations
associated to abelian varieties,
these representations of lifted abelian $D_k^\HH$-type are the only
`new' representations that one may encounter.
To be precise, let $A/F$ be an abelian variety with Mumford--Tate
group $G_A$. 
For any group $G'$ and isogeny
$G^{\prime\der}\rightarrow G_A^\der$ and any weak
geometric lift $\rho'_p\colon\cG_F\rightarrow G'(\Qp)$ of
$\rho_{A,p}$, the representation $\rho'_p$ belongs to the tannakian
category generated by the $p$-adic Galois representations associated
to abelian varieties, the representations 
of lifted abelian $D_k^\HH$-type and those with finite image, see
corollary~\ref{LADHs_dominate_all}. 

In the final section~\ref{sect_LADH} we study the representations $\trho_p$
of lifted abelian $D_k^\HH$-type. 
It is proved in~\ref{LADH_non-abelian} that
such a representation $\trho_p$ generally does not belong to the
tannakian category generated by the Galois representations associated
to abelian varieties. 
The question if $\trho_p$ belongs to the tannakian category
generated by the Galois representations associated to algebraic
varieties, as predicted by the Fontaine--Mazur conjecture, remains open. 

Combined with `standard' conjectures, 
the conjecture of Fontaine and Mazur implies that the $\trho_p$ have
certain `motivic' properties. 
For example, any Fro\-be\-nius element in $\cG_F$ should act semi-simply
on any representation $\wtV_p$ of $\wtG_{/\Qp}$ and the eigenvalues of
the image should be Weil numbers, i.~e.\ algebraic integers of which all
complex absolute values coincide. 
Several of these motivic properties are proved in the final part of
section~\ref{sect_LADH}. 

In addition to their relevance to the Fontaine--Mazur conjecture, the
construction of weak geometric liftings has applications to the
study of Galois representations associated to abelian varieties. 
We have seen that for an abelian variety $A/F$, and for $F$ large
enough, the representations
$\rho_{A,p}\colon\cG_F\rightarrow G_A(\Qp)$ have weak geometric
liftings $\rho_{B,p}\colon\cG_F\rightarrow G_B(\Qp)$, where $B/F$ is
an abelian variety and 
$G_B^\der\rightarrow G_A^\der$
is not far from the universal cover.
The Galois representation on $\rH_\et^1(A_{\bar F},\Qp)$ belongs to
the tannakian category generated by $\rH_\et^1(B_{\bar F},\Qp)$ and
the representations associated to abelian varieties of CM-type. 
This can be useful because the representation 
$\rH_\et^1(B_{\bar F},\Qp)$ is in general easier
to study than 
$\rH_\et^1(A_{\bar F},\Qp)$.
A first example of such an application can be found in
the paper \cite{paugam:good_reduction} of F.~Paugam. 
In \cite{vasiu:M-T_conj}, A.~Vasiu  uses his technique of adjoint
isogenous abelian varieties to study new instances of the
Mumford--Tate conjecture. 
Other applications are to follow in forthcoming publications.
\par\noindent\textbf{Acknowledgements.}
This paper has benefited from discussions and correspondence with a
number of people 
and I thank them heartily for their contribution. 
In particular, I would like to thank D. Blasius and P. Deligne for
(independently) pointing out, after my previous work on the conjecture
of Fontaine and Mazur, that the case of abelian varieties with
Mumford--Tate group of type $D_k^\HH$ deserved special attention. 
\section{Preliminaries}\label{preliminaries}
\subsection{Basic notations.}\label{prel_basic}
For any field $F$, we denote by $\bF$ an algebraic closure of $F$. 
The absolute Galois group of $F$ is the group 
$\cG_F=\Aut_F(\bF)$.

For any prime number $p$, the field $\Qp$ is the $p$-adic completion
of $\QQ$ and $\Cp$ is the completion of an algebraic closure $\bQp$. 

If $G$ is a group and $K$ a field, then $\Rep_K(G)$ is the category of
finite dimensional $K$-linear representations of $G$. 
If $G$ is a topological group (a Galois group for example) and $K$ a
topological field, then the representations in $\Rep_K(G)$ are assumed
to be continuous. 
If $G$ is a linear algebraic group, then $\Rep_K(G)$ is the category
of algebraic representations of $G$. 

If $G$ is an algebraic group, a \emph{quasi-cocharacter} of $G$ is
an element of the direct limit 
\[
\lim_{\genfrac{}{}{0pt}{1}{\longrightarrow}{k\in\NN}}\Hom(\Gm{}^{(k)},G),
\]
where the transition map $\Gm{}^{(k\ell)}\rightarrow\Gm{}^{(k)}$ is the
morphism $z\mapsto z^\ell$. 
Giving a quasi-cocharacter of $G$ is equivalent to giving an integer $k$ and
a cocharacter $\mu\colon\Gm{}\rightarrow G$. 
Intuitively the quasi-cocharacter given by $(k,\mu)$ is the $k$th
root of $\mu$. 
\subsection{Absolute Hodge motives.}\label{prel_AH-mot}
We will freely use the language of tannakian categories. 
Everything we need here can be found in \cite{delmil:tannakian}
which is to be considered the authoritative reference for all notions
used but not explained in this paper. 
In particular, if $\cC$ is a tannakian category, then the subcategory
$\otimes$-generated by a collection $X$ of objects of $\cC$ is the
smallest subcategory of $\cC$ containing all objects which are
isomorphic to a subquotient of a polynomial expression with
coefficients in $\NN$ in the objects
contained in $X$. 
In such a polynomial expression, $+$ and $\cdot$ are to be interpreted
as $\oplus$ and $\otimes$ respectively. 

For any field $F$ of characteristic $0$, 
let $\Mot_\AH(F)$ be the category of motives for absolute
Hodge cycles as described in \cite{delmil:tannakian}, especially
section~6 of that paper. 
It is constructed as Grothendieck's category of Chow motives except where
it concerns the morphisms which are given by absolute Hodge classes,
not by cycle classes as for Grothendieck motives. 
It must be pointed out that by definition, the morphisms between two
motives $M_1$ and $M_2$ defined over $F$ are the appropriate absolute
Hodge classes on the product which are defined over $F$, i.\ e.\ Hodge
classes on $\left(M_1\times M_2\right)_{\bar F}$ invariant under the
action of $\cG_F$. 

Motives are graded objects, each motive $M$ is a finite direct sum
$\bigoplus M^i$, where each $M^i$ is a pure motive of weight $i$. 
The \emph{Tate motive} $\QQ(1)$ is the dual of $h^2(\mathbf{P}^1_F)$ in
the category $\Mot_\AH(F)$. 
\subsection{Mumford--Tate groups.}\label{prelim_MT-group}
In what follows, we will assume that $F$ is contained
in $\CC$ and we write $\bF$ for the algebraic closure of $F$ in
$\CC$. 
Assume that $M_1,\ldots,M_r$ are absolute Hodge motives over $F$, and
let 
$
\cC=\langle M_1,\ldots,M_r,\QQ(1)\rangle
$
be the tannakian subcategory of $\Mot_\AH(F)$ which is
$\otimes$-generated by $M_1,\ldots,M_r$ and the Tate motive $\QQ(1)$. 
Let $\rH_\B(\cdot,\QQ)$ be the fibre functor of $\cC$ over
$\QQ$ defined by the Betti cohomology $\rH^*_\B(X(\CC),\QQ)$ of
complex algebraic varieties.  
Then the \emph{Mumford--Tate group $G$ of $\cC$} is by definition the
automorphism group of $\rH_\B(\cdot,\QQ)$. 
The connected component of $G$ is reductive and $G$ acts on
$\rH_\B(M,\QQ)$ for every object 
$M$ of $\cC$.
The fibre functor $\rH_\B(\cdot,\QQ)$ induces an equivalence between
$\cC$ and the category $\Rep_\QQ(G)$. 
In particular, for $M,M'\in\cC$, one has 
\[
\Hom_{\cC}(M,M')=\Hom_G(\rH_\B(M,\QQ),\rH_\B(M',\QQ)),
\]
i.\ e.\ the action of $G$ fixes all absolute Hodge classes defined
over $F$ on all objects of $\cC$. 
Moreover, $G$ is the smallest $\QQ$-algebraic group with this
property. 
Note that the connected component of $G$ is the Mumford--Tate group of
the subcategory of $\Mot_\AH(\bar F)$ which is $\otimes$-generated by
$M_{1/\bF},\ldots,M_{r/\bF}$ and $\QQ(1)$.
When considering only one motive $M$, the Mumford--Tate group $G_M$
of $M$ is the Mumford--Tate group of the subcategory 
$\langle M,\QQ(1)\rangle$ of $\Mot_\AH(F)$. 

For a general subcategory $\cC$ of $\Mot_\AH(F)$, we define the
Mumford--Tate group in a similar way, obtaining a pro-algebraic
group. 
\subsection{Subcategories of $\Mot_\AH$.}\label{prelim_subcat-mot}
An \emph{Artin motive} is an object of $\Mot_\AH(F)$ with finite
Mumford--Tate group. 
The category of Artin motives is the tannakian subcategory of
$\Mot_\AH(F)$ generated by the finite $F$-schemes. 
An \emph{abelian motive} over $F$ is an AH-motive which belongs to
the tannakian subcategory of $\Mot_\AH(F)$ generated by the
motives of abelian varieties, the Tate motive and the Artin motives. 
We will write $\Artin_F$ for the category of Artin motives and $\AV_F$
for the category of abelian motives over $F$. 
An abelian variety is \emph{potentially of CM-type} if the connected
component of its Mumford--Tate group is commutative and $\CM_F$ is the
tannakian subcategory of $\Mot_\AH(F)$ generated by the Artin motives,
$\QQ(1)$ 
and the motives of the abelian varieties which are potentially of
CM-type. 

If $A$ is an abelian variety over $F$, then 
$h(A)=\bigoplus_i\bigwedge^i h^1(A)$,
so in most questions concerning abelian motives, it suffices to
consider the motives 
$h^1(A)$ instead of the $h(A)$.
\subsection{Betti realization.}\label{prelim_betti}
For any absolute Hodge motive $M$, the Betti realization
$\rH_\B(M,\QQ)$ carries a Hodge structure. 
Giving this Hodge structure is equivalent to giving an action of the
group $S=\CC^\times$, viewed as an algebraic group over $\RR$, on the
$\RR$-vector space $\rH_\B(M,\QQ)\otimes\RR$. 
As the the Mumford--Tate group $G_M$ of $M$ fixes all (absolute) Hodge
classes, this action is given by a morphism of algebraic groups
$h\colon S\rightarrow G_{M/\RR}$.
The couple $(G_M,h)$ is the \emph{Mumford--Tate datum} associated to
$M$. 
More generally, if $\cC$ is a $\otimes$-subcategory of $\Mot_\AH(F)$,
with Mumford--Tate group $G$, then there is a morphism 
$h\colon S\rightarrow G_{/\RR}$ which functorially defines the Hodge
structures on $\rH_\B(M,\QQ)$ for all objects $M$ of $\cC$. 

We restrict our attention to the category of abelian motives $\AV_F$. 
As explained in \cite[6.25]{delmil:tannakian}, it follows from
\cite[Theorem~2.11]{deligne:hodge_cycles} that the Betti
realization induces an equivalence of $\AV_F$ with its essential image
in the category $\Hodge$ of Hodge structures. 

Let $\cC$ be a $\otimes$-subcategory of $\Mot_\AH(F)$ and let
$(G,h)$ be the associated Mumford--Tate datum. 
For any object $M$ of $\cC$, the Betti realization $\rH_\B(M,\QQ)$ is a
representation of $G$ and the Hodge structure on $\rH_\B(M,\QQ)$ is
given by the action of $S$ on $\rH_\B(M,\QQ)\otimes\RR$ induced by
$h\colon S\rightarrow G_{/\RR}$. 
Since every $\QQ$-linear representation of $G$ is the Betti
realization of an object of $\cC$, 
this construction defines a $\otimes$-equivalence of $\Rep_\QQ(G)$
with the essential image of $\cC$ in $\Hodge$. 

If $A$ is an abelian variety, its Mumford--Tate group $G_A$ is equal
to the Mum\-ford--Tate group of $h^1(A)$.
The connected component of $G_A$ coincides in turn with
the Mumford--Tate group of the Hodge structure $\rH_\B^1(A(\CC),\QQ)$. 

Let $(G,h)$ be a Mumford--Tate datum. 
Composition of $h$ with the cocharacter 
$\Gm{/\CC}\rightarrow S_{/\CC}$ dual to 
$z\colon S_{/\CC}\rightarrow\Gm{/\CC}$
gives rise to \emph{Hodge cocharacter}
$\mu\colon\Gm{/\CC}\rightarrow G_{M/\CC}$. 
Alternatively, $\mu$ is determined by the condition that $\Gm{/\CC}$
acts on the factor $V^{p,q}$ of the Hodge decomposition as
multiplication by $z^p$. 
Conversely, a pure Hodge structure is determined by giving its weight
and the Hodge cocharacter. 
\subsection{$p$-adic Galois representations.}\label{prelim_p-adic}
Let $p$ be a prime number. 
The \'etale cohomology with coefficients in $\Qp$ of algebraic
varieties over $\bar F$ defines the $p$-adic \'etale realization on
the category of absolute Hodge motives. 
The $p$-adic \'etale realization of a motive $M$ over $F$ is a
$\Qp$-vector space $\rH_\et(M_{\bar F},\Qp)$ endowed with a continuous
action of the group $\cG_F$.
It follows from standard tannakian theory that the image of $\cG_F$ in
$\GL(\rH_\et(M_{\bar F},\Qp))$ lies in $G_M(\Qp)$, where $G_M$ is the
Mumford--Tate group of $M$. 
The representation thus gives rise to a morphism
$\rho_{M,p}\colon\cG_F\rightarrow G_M(\Qp)$.

For any prime number $p$, let $\Artin\dash\Rep_{\Qp}(\cG_F)$ be the 
tannakian subcategory of $\Rep_{\Qp}(\cG_F)$ generated by the
$p$-adic \'etale realizations of the objects of $\Artin_F$. 
This category coincides with the category of $\Qp$-linear
representations of $\cG_F$ with finite image. 
Similarly, we let let $\RepAV{\Qp}{(\cG_F)}$ be the 
tannakian subcategory of $\Rep_{\Qp}(\cG_F)$ generated by the
$p$-adic \'etale realizations of the objects of $\AV_F$.
As before, in the case of an abelian variety $A$, 
it is usually sufficient to consider just
the representation on the first \'etale cohomology group
$\rH_\et^1(A_{\bar F},\Qp)$. 

Restricting to an inertia group at a $p$-adic place of $F$, these
representations give rise to representations of 
\emph{Hodge--Tate type}, cf.~\cite[\S3]{fontaine:repr_semist} or
\cite{serre:hodge_tate}. 
If $\cI\subset\cG_F$ is such an inertia subgroup and $V_p$ is a
$\Qp$-linear representation of $\cI$ of Hodge--Tate type, there is a
canonical decomposition
\[
V_p\otimes_{\Qp}\Cp=\bigoplus_{p,q}V^{p,q}, 
\]
the \emph{Hodge--Tate decomposition} of $V_p\otimes\Cp$.
The \emph{Hodge--Tate cocharacter} is the cocharacter
$\mu\colon\Gm{/\Cp}\rightarrow\GL(V_p\otimes\Cp)$ such that action of
$\Gm{/\Cp}$ on $V^{p,q}$ is the multiplication by $x^p$. 
It can be shown (see \cite[\S1]{serre:hodge_tate}) that the connected
component $H$ of the Zariski closure of the image of $\cI$ in
$\GL(V_p)$ coincides with the smallest subgroup $H\subset\GL(V_p)$
such that $\mu$ factors through $H_{/\Cp}$. 
\section{Mumford--Tate liftings of abelian varieties over $\CC$}
\label{sect_MT_liftings}
\subsection{Mumford--Tate liftings.}\label{def_MT-lift}
Let $M/\CC$ be an abelian motive, $G_M$ its Mumford--Tate group
and $(G_M,h_M)$ the associated Mumford--Tate datum. 
This means that the Hodge structure on $V_M=\rH_\B(M,\QQ)$ is defined by  
$h_M\colon S\rightarrow G_{M/\RR}$. 

We say that an abelian motive $N$ with Mumford--Tate datum
$(G_N,h_N)$ pro\-vides a \emph{Mumford--Tate lift} of $M$ if there
exists a central isogeny 
$\pi\colon G_N\rightarrow G_M$ such that 
$\pi_\RR\circ h_N=h_M$.
We say that $M$ is \emph{Mumford--Tate liftable} if there exists an
abelian motive $N/\CC$ giving a Mumford--Tate lift of $M$ and such
the morphism $\pi\colon G_N\rightarrow G_M$ is not an isomorphism.  
We say that $M$ is \emph{Mumford--Tate unliftable} if it is not
Mumford--Tate liftable. 

If there exists a central isogeny 
$\pi^\der\colon G_N^\der\rightarrow G_M^\der$ such that 
\[\pi_N^\ad\circ h_N=\pi_M^\ad\circ h_M\] then 
$N$ provides a \emph{weak Mumford--Tate lift} of $M$.
Here the maps $\pi_M^\ad$ and $\pi_N^\ad$ are the projections
$G_M\rightarrow G_M^\ad$ and $G_N\rightarrow G_N^\ad$. 
As it is a central isogeny, it follows that $\pi^\der$ induces an isomorphism 
$G_M^\ad\cong G_N^\ad$, giving a sense to the equality 
$\pi_N^\ad\circ h_N=\pi_M^\ad\circ h_M$.
Finally, $M$ is \emph{essentially Mumford--Tate unliftable} if there
does not exists any abelian motive $N/\CC$ giving a weak
Mumford--Tate lift of $M$ for which $\pi^\der$ is not an isomorphism. 

We will often write M-T (un)liftable instead of
Mumford--Tate (un)liftable. 
\subsection{Remark.}
If $N$ is a weak Mumford--Tate lift of $M$, then $M$ and $N$ are
adjoint-isogenous in the sense of \cite[6.1]{vasiu:M-T_conj} and
conversely, if $M$ and $N$ are adjoint-isogenous and if there is an
isogeny $G_N^\der\rightarrow G_M^\der$, then $N$ is a weak
Mumford--Tate lift of $M$. 
\subsection{The Mumford--Tate datum of an abelian variety.}
\label{structure_M-T}
Let $A/\CC$ be an abelian variety, 
$(G_A,h)$ the associated Mum\-ford--\-Tate datum and 
$\mu\colon\Gm{/\CC}\rightarrow G_{A/\CC}$ be the Hodge cocharacter. 
It follows from \cite{deligne:shimura2}, in particular from~1.3
and~2.3, that the simple factors of $G^\ad_{A/\RR}$ are absolutely simple
and of classical type ($A$, $B$, $C$ or $D$). 
The group $\cG_\QQ$ acts on the set of these simple factors 
and each factor is conjugate to a non-compact one. 
The compact factors are exactly the factors to which $h$ projects
trivially. 

Let $H$ be a $\QQ$-simple factor $G_A^\ad$ and decompose 
$H_{/\CC}=\prod_{\iota\in I_H}H_{/\CC,i}$ for some finite set $I_H$
with $\cG_\QQ$-action. 
Here we have identified $\bQ$ with the algebraic closure of $\QQ$ in
$\CC$. 
Note that the complex conjugation in $\cG_\QQ$ acts trivially on $I$
because the simple factors of $G^\ad_{A/\RR}$ are absolutely simple.
Let $I_{H,c}$ be the set of indices such that the Hodge cocharacter
$\mu$ projects trivially to $H_{/\CC,\iota}$ and let $I_{H,nc}$ be the
complement of $I_{H,c}$ in $I_H$. 
The $\iota\in I_{H,c}$ are exactly the indices for which the
corresponding real factor of $H_{/\RR}$ is compact. 
For each $\iota\in I_{H,nc}$, the Hodge cocharacter lifts to a
quasi-cocharacter $\tmu_\iota$ of the universal cover
$\wtH_{/\CC,\iota}$ of $H_{/\CC,\iota}$.  
If the simple factors of $H_{/\CC}$ are of type $A$, $B$ or $C$ then
it follows from \cite[1.3]{deligne:shimura2} that $\wtH_{/\CC,\iota}$
admits a faithful representation $W_{/\CC,\iota}$ such that
$\tmu_\iota$ acts on $W_{/\CC,\iota}$ with exactly two weights. 
Contemplating the tables \cite[1.3.9]{deligne:shimura2} or
\cite[Table~4.2]{pink:monodromy_groups} one sees that
the highest weight of $W_{/\CC,\iota}$ is
\begin{itemize}
\item
either $\varpi_1$ or $\varpi_k$ if the simple factors of $H_{/\CC}$
are of type $A_k$, 
\item
$\varpi_k$ if these simple factors are of type $B_k$ and 
\item
$\varpi_1$ if the factors are of type $C_k$.
\end{itemize}
For each $\iota\in I_H$, we define a representation $V_{/\CC,\iota}$
of $\wtH_{/\CC,\iota}$ as follows.
Let $V_{/\CC,\iota}$ be the direct sum of the representations with
highest weights $\varpi_1$ and $\varpi_k$ if $H$ is of type $A_k$ with
$k\geq2$ and define $V_{/\CC,\iota}$ to be the representation with
highest weight $\varpi_1$ (resp.\ $\varpi_k$) if $H$ is of type $A_1$
or $C_k$  (resp.\ $B_k$).
For $\iota\in I_{H,nc}$, the cocharacter $\tmu_\iota$ still acts on
each irreducible factor of $V_{/\CC,\iota}$ with exactly two weights. 
The product over $\iota\in I_H$ of the $\wtH_{/\CC,\iota}$
descends to an algebraic group $\wtH$ over $\QQ$ and a multiple of 
the direct sum of the $V_{/\CC,\iota}$ descends to a faithful
$\QQ$-linear representation $V$ of $\wtH$. 

\subsection{Remarks.}\label{remarks_V_iota}
If the simple factors of $H_{/\CC}$ are of type $B_k$ or $C_k$ then
the condition that $\tmu_\iota$ acts on $W_{/\CC,\iota}$ with exactly
two weights uniquely determines the highest weight.
Similarly, if the simple factors of $H_{/\CC}$ are of type $A_k$ and
if $\tmu_\iota$ is not dual to $\alpha_1$, only the representations
with highest weight $\varpi_1$ or $\varpi_k$ fulfill this
condition. 
Obviously, all these representations are faithful.

In the case where the simple factors of $H_{/\CC}$ are of type $A_k$
and where $\tmu_\iota$ is dual to $\alpha_1$, 
the cocharacter $\tmu_\iota$ acts with exactly two weights on 
the representation with highest weight $\varpi_s$ for any 
$1\leq s\leq k$, see \cite[Table~4.2]{pink:monodromy_groups}.
Only for $s=1$ and for $s=k$ does one obtain a faithful
representation. 

Also note that $V_{/\CC,\iota}$ is self dual in each case, see the
aforementioned table in Pink's paper.

\subsection{Factors of type $D_k$.}\label{D_k-factors}
In what follows, we will focus on the factors of $G_A^\ad$ of type $D_k$, first
considering the case $k\geq5$. 
Fix a $\QQ$-simple factor $H$ of $G_A^\ad$ such that the simple factors
of $H_{/\CC}$ are of type $D_k$ with $k\geq5$. 
Let $H_{/\CC}=\prod_{\iota\in I_H}H_{/\CC,\iota}$ and let $I_{H,c}$ and
$I_{H,nc}$ be as before. 
The Dynkin diagram of $H_{/\CC}$ is the disjoint union, indexed by $I_H$, of
diagrams of type $D_k$, and it follows from
\cite[1.3]{deligne:shimura2} that, for each $\iota\in I_{H,nc}$, the
conjugacy class of the cocharacter $\mu_\iota$ is dual to an endpoint
$\alpha^{(\iota)}$ of the $\iota$-component of the Dynkin diagram. 
\begin{lemma}
Assume that $k\geq5$ and that
$\alpha^{(\iota)}=\alpha_1$ for some $\iota\in I_{H,nc}$. 
Then
$\alpha^{(\kappa)}=\alpha_1$ for all $\kappa\in I_{H,nc}$. 
\end{lemma}
\begin{proof}
Assume that $\alpha^{(\iota)}=\alpha_1$ for $\iota\in I_{H,nc}$ and let
$G^\der_{A/\CC,\iota}$ be the simple factor of $G^\der_{A/\CC}$ projecting onto
$H_{/\CC,\iota}$. 
Via the representation of $G_A$ on $V=\rH_\B^1(A(\CC),\QQ)$, the Hodge
cocharacter $\mu$ acts on $V_{/\CC}$ with two weights. 
By the table in \cite[1.3]{deligne:shimura2}, this implies that every
non-trivial irreducible direct factor of the 
representation of $G^\der_{A/\CC,\iota}$ on $V_{/\CC}$ has highest weight
$\varpi_{k-1}$ or $\varpi_k$. 
As $H$ is $\QQ$-simple, the same thing is true for the
representation on $V_{/\CC}$ of the simple factors of $G^\der_{A/\CC}$
mapping to the factors $H_{/\CC,\kappa}$ for the other $\kappa\in I_H$. 
Using the tables \cite[1.3.9]{deligne:shimura2},
the fact that $\mu$ acts with two weights implies that
all non-trivial $\mu_\kappa$ are dual to $\alpha_1$. 
\end{proof}
If any factor of $H_{/\CC}$ satisfies the conditions of the lemma, then
we say that $H$ and the factors of $H_{/\RR}$ and $H_{/\CC}$ are
\emph{of type $D_k^\RR$}.
As $H$ was assumed to be $\QQ$-simple, the vertices
$\alpha_1^{(\iota)}$ form a single orbit for the action of $\cG_{\bQ}$
in this case. 
In the opposite case we say that they are \emph{of type $D_k^\HH$}. 
In the latter case, the conjugacy class of the projection of $\mu$ to
any factor of $H_{/\CC}$ is either trivial or dual to $\alpha_{k-1}$
or $\alpha_k$. 

Let $H_{/\CC,\iota}$ be a factor of $G^\ad_{A/\CC}$ of type $D_k^\RR$
and let $\wtH_{/\CC,\iota}$ be its universal cover.
For each $\iota\in I_H$, let $V_{/\CC,\iota}$ be the direct sum of the
representations of $H_{/\CC,\iota}$ with highest weights
$\varpi_{k-1}$ and $\varpi_k$. 
As above, for $\iota\in I_{H,nc}$, the projection of the Hodge
cocharacter $\mu$ to $H_{/\CC,\iota}$ lifts to a quasi-cocharacter
$\tmu_\iota$ of $\wtH_{/\CC,\iota}$ and $\tmu_\iota$ acts on
$V_{/\CC,\iota}$ with two weights $\pm1/2$. 
The product over $\iota\in I$ of the $\wtH_{/\CC,\iota}$ and a multiple
of the direct sum of the $V_{/\CC,\iota}$ descend to an algebraic
group $\wtH$ over $\QQ$ and a faithful $\QQ$-linear representation $V$ of
$\wtH$. 

Next assume that $H_{/\CC,\iota}$ is a factor of $G^\ad_{A/\CC}$ of
type $D_k^\HH$. 
As before, for $\iota\in I_{H,nc}$, the projection of the Hodge
cocharacter lifts to a quasi-cocharacter $\tmu_\iota$ of the universal
cover $\wtH_{/\CC,\iota}$, but this time the group $\wtH_{/\CC,\iota}$
does not have a faithful representation on which $\tmu_\iota$ acts
with two weights. 
This property can only be achieved for a quotient of
$\wtH_{/\CC,\iota}$ which can be constructed as follows. 
For each $\iota\in I_H$, let $V_{/\CC,\iota}$ be the representation of
$\wtH_{/\CC,\iota}$ with highest weight $\varpi_1$. 
We will refer to the quotient of $\wtH_{/\CC,\iota}$ which acts
faithfully on $V_{/\CC,\iota}$ as the \emph{$h$-maximal cover} of
$H_{/\CC,\iota}$.  
This $h$-maximal cover descends over $\RR$ to an inner form of the
orthogonal group $\SO(2k)$.
For each $\iota\in I_{H,nc}$, the quasi-cocharacter $\tmu_\iota$ acts
on $V_{/\CC,\iota}$ with weights $\pm1/2$. 
The product over $\iota\in I$ of these $h$-maximal covers descends to
an algebraic group over $\QQ$ and a multiple of the direct sum of
the $V_{/\CC,\iota}$ descends to a faithful $\QQ$-linear
representation of this group. 
\subsection{Factors of type $D_4$.}\label{D_4-factors}
We finally consider the case where $k=4$, so let $H$ be a $\QQ$-simple
factor of $G_A^\ad$ such that the simple factors 
of $H_{/\CC}$ are of type $D_4$. 
In this case, the automorphism group of the Dynkin diagram permutes
the set of endpoints, so the types $D_4^\RR$ and $D_4^\HH$ can not be
distinguished in the same manner as before. 
The fundamental difference between the two types is the existence of a
$\QQ$-linear representation $V$ as above. 
We formalize this as follows. 

Let $H_{/\CC}=\prod_{\iota\in I_H}H_{/\CC,\iota}$ and 
$I_{H,nc}$ and $I_{H,c}$ as in the general case. 
For each $\iota\in I_{H,nc}$ let 
$\mu_\iota\colon\Gm{/\CC}\rightarrow H_{/\CC,\iota}$ be the projection
of the Hodge cocharacter. 
The conjugacy class of each $\mu_\iota$ is dual to an
endpoint of the corresponding component of the Dynkin diagram. 
Let $\Delta$ be the set of these vertices, it contains exactly one 
endpoint of the $\iota$-component of the Dynkin diagram if 
$\iota\in I_{H,nc}$ and no vertices in the other components.
Choose $\iota\in I_{H,nc}$, let $G^\der_{A/\CC,\iota}$ be the almost
simple factor of 
$G^\der_{A/\CC}$ projecting onto $H_{/\CC,\iota}$ and let $\varpi$ be the
highest weight of some non trivial factor of the representation of
$G^\der_{A/\CC,\iota}$ on $\rH_\B^1(A(\CC),\QQ)$. 
As the Hodge cocharacter acts on $\rH_\B^1(A(\CC),\CC)$ with two
weights, the arguments of \cite[1.3, 2.3]{deligne:shimura2} imply that
$\varpi$ corresponds to an endpoint $\alpha$ of the Dynkin diagram of
$H$ and that the $\cG_\QQ$-orbit $\Delta'$ of $\alpha$ does not meet
$\Delta$. 
In particular, there exists a $\cG_\QQ$-stable set $\Delta'$ of
endpoints such that $\Delta\cap\Delta'=\emptyset$.
The set $\Delta'$ contains at least one endpoint in each connected
component of the Dynkin diagram. 

If there exists a $\cG_\QQ$-stable set $\Delta'_{\max}$ of endpoints of the
Dynkin diagram with $\Delta\cap\Delta'_{\max}=\emptyset$ containing two
endpoints in each connected component, then we say that
$H$ and the factors of $H_{/\RR}$ and $H_{/\CC}$ are 
\emph{of type $D_4^\RR$}.
If there does not exist such a $\Delta'_{\max}$, they are 
\emph{of type $D_4^\HH$}.  

Assume that $H$ is a factor of $G^\ad_{A/\CC}$ of type $D_4^\RR$, and
let $\iota\in I_H$. 
Then the $\iota$-component of the Dynkin
diagram contains two endpoints $\alpha_\iota$ and
$\beta_\iota\in\Delta'_{\max}$. 
Let $V_{/\CC,\iota}$ be the direct sum of the representations of the
universal cover $\wtH_{/\CC,\iota}$ with highest weights
corresponding to $\alpha_\iota$ and $\beta_\iota$ respectively. 
It is a faithful representation of $\wtH_{/\CC,\iota}$ and as
$\alpha_\iota,\beta_\iota\not\in\Delta$, it follows from
\cite[1.3.9]{deligne:shimura2} that, for $\iota\in I_{H,nc}$, the
lifting $\tmu_\iota$ of $\mu_\iota$ to $\wtH_{/\CC,\iota}$ acts on
$V_{/\CC,\iota}$ with weights $\pm1/2$. 
The product of the $\wtH_{/\CC,\iota}$ and a multiple of the direct
sum of the $V_{/\CC,\iota}$ descend to a $\QQ$-algebraic group $\wtH$
with a faithful $\QQ$-linear representation $V$. 

If $H_{/\CC,\iota}$ is a factor of $G^\ad_{A/\CC}$ of type $D_4^\HH$,
then the above set $\Delta'$ meets the Dynkin
diagram of $H_{/\CC,\iota}$
in one endpoint, corresponding to a fundamental weight $\varpi$. 
The \emph{$h$-maximal cover} of $H_{/\CC,\iota}$ is the quotient of its
universal cover $\wtH_{/\CC,\iota}$ acting faithfully on the
representation 
of $\wtH_{/\CC,\iota}$ with highest weight $\varpi$. 
As in the case $k\geq5$, the product of these data over all $\iota\in
I_H$ descends to an
algebraic group over $\QQ$ together with a faithful $\QQ$-linear
representation. 
This cover of a $\QQ$-simple factor of $G_A^\ad$ of type $D_4^\HH$ is
the $h$-maximal cover of $G_A^\ad$.
\subsection{Remark.}
For all factors $H_{/\CC,\iota}$ of type $D$, the representation
$V_{/\CC,\iota}$ is self dual. 
As in the case of the factors of types $A$, $B$ and $C$ this can be
read off from \cite[Table~4.2]{pink:monodromy_groups}.
\begin{thm}\label{MT-group_unliftable-AV}
Let $A/\CC$ be an abelian variety and
let $(G_A,h)$ be the associated Mum\-ford--\-Tate datum. 
Then the following conditions are equivalent.
\begin{enumerate}
\item\label{prop12-point1}
$A$ is essentially Mumford--Tate unliftable.
\item\label{prop12-point2}
The group $G_{A/\CC}^\der$ is the product of its simple factors.
The simple factors of types $A_k$, $B_k$, $C_k$ and $D_k^\RR$ are
simply connected and the factors of type $D_k^\HH$ are $h$-maximal in
the sense defined above. 
\end{enumerate}
\end{thm}
\begin{dfn}\label{def_MT-reduced}
Let $A/\CC$ be an abelian variety, $G_A$ its
Mumford--Tate group and $V_A=\rH_\B^1(A(\CC),\QQ)$. 
If $G_A^\ad$ is $\QQ$-simple, then we say
that $A$ is \emph{Mumford--Tate decomposed} if the following 
conditions hold. 
\begin{itemize}
\item
There are a totally real field $K_0$ and an absolutely simple
algebraic group $G^s$ over $K_0$ such that
$G_A^\der\cong\Res_{K_0/\QQ}G^s$. 
\item
There is a faithful representation $V^s$ of $G^s$ such that the
representation of $G_A^\der$ on $V$ is isomorphic to
$\Res_{K_0/\QQ}V^s$. 
\item
There is no proper non-trivial abelian subvariety $B$ of $A$ with
Mumford--Tate group $G_B$ verifying $G_B^\der=G_A^\der$. 
\end{itemize}

We say that an abelian variety $A/\CC$ 
is \emph{Mumford--Tate decomposed} if it is isogenous to a product 
$\prod A_i$ such that 
each $A_i$ is Mumford--Tate decomposed with $G_{A_i}^\ad$ simple
and if $G_A^\ad$ is the
product of the $G_{A_i}^\ad$. 
An abelian variety over a number field $F\subset\CC$ is
\emph{Mumford--Tate decomposed} if $A_{/\CC}$ is Mumford--Tate decomposed.
\end{dfn}
\subsection{Remark.}
The upshot of this definition is that the abelian varieties arising
from the construction of \cite[2.3]{deligne:shimura2} are
Mumford--Tate decomposed. 
The representation of $G_A^\der$ on $\rH_\B^1(A(\CC),\QQ)$ is a
direct sum of the representations of the $G_{A_i}^\der$ constructed
in~\ref{structure_M-T} and~\ref{D_k-factors}. 
See also the proofs given below.

The notion of being (essentially) Mumford--Tate unliftable
is a condition 
on the Mumford--Tate datum of an abelian variety whereas the notion of
being Mumford--Tate decomposed pertains to the action of the
Mumford--Tate group on the first Betti cohomology group. 
\begin{thm}\label{unliftable_lift_exists}
For every abelian variety $A/\CC$ there exists a weak
Mumford--Tate lift $B/\CC$ of $A$ such that 
$B$ is essentially Mumford--Tate unliftable and 
Mumford--Tate decomposed.
\end{thm}
\begin{proof}[Proofs of~\ref{MT-group_unliftable-AV}
and~\ref{unliftable_lift_exists}] 
These results can be derived from the work of Satake,
see for example~\cite{deligne:shimura2}. 
The same argument can be found in \cite{vasiu:M-T_conj},
see \S4 and paragraphs~6.3 and~6.4 in particular.
The strategy of the proof is as follows. 
It is first shown that the condition \ref{prop12-point2} implies the
condition \ref{prop12-point1}. 
We then prove that any abelian variety admits a weak Mumford--Tate lift
satisfying \ref{prop12-point2} and which is Mumford--Tate decomposed. 
Thanks to the fact that \ref{prop12-point2} implies
\ref{prop12-point1}, this M-T lift is also M-T unliftable. 
Finally, the proof that the condition \ref{prop12-point1} implies
\ref{prop12-point2} is a formality.

First assume that $A$ verifies the conditions of \ref{prop12-point2}. 
We will show that it is essentially M-T unliftable, so
let $B$ be a weak Mumford--Tate lift of $A$ and let $G_B$ be its 
Mumford--Tate group. 
It has to be proved that the map $G_B^\der\rightarrow G_A^\der$ is an
isomorphism. 
It suffices to prove this after extension of scalars to $\CC$ and as
the only non-simply connected factors of $G_{A/\CC}^\der$ are the
factors of type $D_k^\HH$, we only need to consider these factors. 

Consider a factor $H_{/\CC}$ of $G^\der_{A/\CC}$ of type $D_k^\HH$
to which the Hodge cocharacter projects non-trivially, 
assuming at first that $k\geq5$. 
As we saw, the conjugacy class of the projection of the Hodge
cocharacter to $H_{/\CC}$ is dual to one of  the vertices
$\alpha_{k-1}$ or $\alpha_k$ of the Dynkin diagram.
Let $\wtH_{/\CC}$ be the factor of $G_{B/\CC}^\der$ mapping onto
$H_{/\CC}$. 
An appropriate direct factor $W$ of the representation of $\wtH_{/\CC}$ on
$\rH_\B^1(B(\CC),\QQ)\otimes\CC$ provides a faithful representation on
which $\tmu$ acts with two weights $\pm1/2$. 
As explained in~\ref{D_k-factors}, 
it follows from the tables in~\cite[1.3.9]{deligne:shimura2}, that
the highest weight of every irreducible direct factor of $W$ is
$\varpi_1$ and hence that 
$\wtH_{/\CC}$ is isomorphic to $H_{/\CC}$. 
This proves that~\ref{prop12-point2} implies~\ref{prop12-point1} if
$k\geq5$.

In the case where $k=4$, first note that $G_A^\der$ is the product of
its $\QQ$-simple factors. 
Let $H$ be a simple factor of type $D_4^\HH$ and let $\wtH$ be the
factor of $G_B$ mapping onto $H$. 
Consider $\rH_\B^1(B(\CC),\QQ)$ as a representation of $\wtH$ and let 
$W$ be a direct factor which is a faithful representation of $\wtH$. 
Let $\wtH_{/\CC,\iota}$ be any factor of $\wtH_{/\CC}$ and let
$W_{/\CC,\iota}$ be any irreducible direct factor of the restriction of
$W\otimes\CC$ to $\wtH_{/\CC,\iota}$.
Then the lifting to $\wtH_{/\CC,\iota}$ of the Hodge cocharacter
either acts trivially on $W_{/\CC,\iota}$ or with exactly two
weights. 
It follows from~\ref{D_4-factors} that $\wtH$ is the $h$-maximal cover
of $H$ and as $H$ was $h$-maximal by hypothesis it follows that
$\wtH\cong H$. 
This proves that ~\ref{prop12-point2} implies~\ref{prop12-point1} in
case $k=4$.

We next show that if $A$ is any abelian variety, then there exists a
weak Mum\-ford--Tate lift $B$ of $A$ with 
Mumford--Tate group $G_B$ satisfying the conditions of
\ref{prop12-point2} and which is Mumford--Tate decomposed. 
This fact readily follows from~\cite[2.3]{deligne:shimura2} and the
discussions in~\ref{structure_M-T}, \ref{D_k-factors}
and~\ref{D_4-factors}, we recall the argument. 

Fix a $\QQ$-simple factor $H^\ad$ of $G_A^\ad$.
It is of the form
$\Res_{K_0/\QQ}H^{s,\ad}$ for some totally real number field $K_0$
and some absolutely simple adjoint group $H^{s,\ad}$ over $K_0$.
As usual, 
decompose $H_{/\CC}^\ad=\prod_{\iota\in I_H} H_{/\CC,\iota}$, where
$I_H$, $I_{H,nc}$ and $I_{H,c}$ are as before. 

Unless $H^\ad_{/\CC,\iota}$ is of type $D_k^\HH$, we let
$\wtH^{s,\der}$ be the simply connected cover of $H^{s,\ad}$. 
In the remaining case we let $\wtH^{s,\der}$ be the $h$-maximal cover
of $H^{s,\ad}$, in the sense of~\ref{D_k-factors}
resp.~\ref{D_4-factors}. 
In all cases, put \[\wtH^\der=\Res_{K_0/\QQ}\wtH^{s,\der}.\]
For each $\iota\in I_{H,nc}$, the projection of the Hodge cocharacter
to $H^\ad_{/\CC,\iota}$ lifts to a quasi-cocharacter $\tmu_\iota$  
of $\wtH^\der_{/\CC,\iota}$.
The faithful representation $V$ of $\wtH^\der$ constructed
in~\ref{structure_M-T} and~\ref{D_k-factors} resp.~\ref{D_4-factors}
is of the form $V=\Res_{K_0/\QQ}V^s$ for some faithful $\QQ$-linear
representation $V^s$ of $\wtH^{s,\der}$. 
These data have the following
properties. 
\begin{itemize}
\item
There is a map $\wtH^\der\rightarrow G_A^\der$ such that
$\tmu=\prod_{\iota\in I_{H,nc}}\tmu_\iota$ lifts the Hodge
cocharacter. 
\item
If $\iota\in I_{H,nc}$, and if $W$ is an irreducible factor of
$V_{\CC,\iota}$, then the cocharacter $\tmu$ acts either trivially on
$W$ or with two rational weights $r$ and $r+1$.
\end{itemize}
There exist a  torus $T$ over $\QQ$ acting $\wtH^\der$-linearly on
$V$ and a quasi-cocharacter $\mu_T$ of $T_{/\CC}$ such that the product
$\tmu\mu_T$ acts trivially on the $V_{\CC,\iota}$ for 
$\iota\in I_{H,c}$ and with weights $\pm1/2$ for $\iota\in I_{H,nc}$. 
In fact, $T$ is characterized by the condition that $T_{/\CC}$ is
the group of the automorphisms of $V\otimes\CC$ acting
by scalar multiplication on each isobaric component of $V\otimes\CC$
as $\wtH^\der_{/\CC}$-representation. 
The existence of $\mu_T$ follows from the fact that if $\tmu$ acts
non-trivially on an isobaric component then it acts with two weights $r$
and $r+1$. 
Let $H'$ be the image of $\wtH^\der\times T$ in $\GL(V)$ and let 
$\mu'=\tmu\mu_T$.

We choose a quadratic and totally imaginary extension $L$ of $K_0$ and
consider the $\QQ$-algebraic group $L^\times$. 
The natural action of $L^\times$ on $L$ gives rise to a $\QQ$-linear
representation $W$ of $L^\times$. 
One has 
\[
\bigl(L^\times\bigr)_{/\CC}\cong\bigoplus_{\iota\in I}\Gm{/\CC}^2
\]
and one defines a quasi-cocharacter $\nu$ of 
$\bigl(L^\times\bigr)_{/\CC}$ by $\nu_\iota(z)=(z^{1/2},z^{1/2})$ for
$\iota\in I_{nc}$ and $\nu_\iota(z)=(z,1)$ for $\iota\in I_c$. 
The action of $H'\times L^\times$ on
$V\otimes_{K_0}W$ defines a faithful representation of a
quotient $\wtH$ of $H'\times L^\times$ with derived group
$\wtH^\der$ in which $(\tmu,\nu)$ acts with weights $0$ and $1$. 
Let $\tih\colon S\rightarrow\wtH_\RR$ be defined by 
\[
\tih(z,\bar z)=(\mu'\nu)(z)\overline{(\mu'\nu)(\bar z)}.
\] 
Shrinking the centre of $\wtH$, we may assume that the image of $\tih$
is Zariski dense in $\wtH$. 
We claim that this defines the Shimura datum $(\wtH,\tih)$ associated
to an abelian variety $B_1$. 
To see why this is the case, note that $\tih$ defines a $\QQ$-Hodge
structure of type $(1,0),(0,1)$ on $V\otimes_{K_0}W$. 
In order to establish that this Hodge structure comes from an abelian
variety it is sufficient to show that it is polarizable, 
cf.\ \cite[2.3]{deligne:K3}. 
First use that, by loc.\ cit.~2.11, the element $\ad(h_A(i))$ defines
a Cartan involution of $G_A^\der$ and hence 
that $\ad(\tih(i))$ is a Cartan involution of
$\wtH_{/\RR}^\der$. 
Next, one checks that the weight $w_h\colon\Gm{}\rightarrow\wtH$ is
central, defined over $\QQ$ and that the quotient of the center of
$\wtH$ by $w_h(\Gm{})$ is compact. 
The last statement is deduced from the fact that this center is
contained in a product of CM tori. 
It now follows that $\ad(\tih(i))$ is a Cartan involution of
$\wtH/w(\Gm{})$ and 
by \cite[1.1.18(b)]{deligne:shimura2} this implies that the Hodge
structure on $V\otimes_{K_0}W$ is polarizable.

The group $\wtH^\ad$ is $\QQ$-simple, $\wtH^\der$ is its $h$-maximal
cover and $B_1$ is Mumford--Tate decomposed by construction. 
Since $G_A$ is the Mumford--Tate group of an abelian variety,
$G_A^\der$ is a quotient of the $h$-maximal cover of $G_A^\ad$, so
there is an isogeny $\wtH^\der\rightarrow G_A^\der$ lifting
$\wtH^\ad\rightarrow G_A^\ad$. 

Applying this to all $\QQ$-simple factors of $G_A^\ad$ we obtain
Mumford--Tate decomposed 
abelian varieties $B_i$ such that $B=\prod B_i$ is a weak M-T lift of $A$
verifying the conditions of~\ref{prop12-point2}. 
By construction, if $G_B$ (resp.\ $G_{B_i}$) denotes the Mumford--Tate
group of $B$ (resp.\ $B_i$), then $G_B^\der=\prod G_{B_i}^\der$.
The first part of this proof implies that $B$ is essentially M-T
unliftable. 
This terminates the proof of theorem~\ref{unliftable_lift_exists}.

Finally, let $A$ be essentially M-T unliftable, i.\ e.\ 
the condition~\ref{prop12-point1} is satisfied. 
The theorem~\ref{unliftable_lift_exists} implies that $A$ has a
weak M-T lift $B$ with Mumford--Tate group $G_B$ satisfying the
condition of~\ref{prop12-point2}. 
As $A$ is essentially M-T unliftable, we must have 
$G_A^\der\cong G_B^\der$, which implies that $G_A$ also
verifies~\ref{prop12-point2}. 
\end{proof}
\subsection{Remarks.}
\subsubsection{}
In the above proof, fix a $\QQ$-simple factor $H^\ad$ of $G_A^\ad$, let
$\iota\in I_H$ and consider an irreducible 
factor $W$ of the representation $V_{\CC,\iota}$ of
$\wtH_{/\CC,\iota}^\der$. 
It follows from \cite[1.3]{deligne:shimura2} that the highest weight
of $W$ is a fundamental weight of $\wtH_{/\CC,\iota}^\der$. 
More precisely, according to the type of $\wtH_{/\CC,\iota}^\der$ and
the quasi-cocharacter $\tmu_\iota$, it
is the weight given by the tables~1.3.9 of Deligne's paper or
\cite[Table~4.2]{pink:monodromy_groups}. 
\subsubsection{}\label{better_unliftable_lift}
With the same notations, assume that $H^\ad$ is not of type $A_k$ with
$k\geq2$. 
In the above construction of the essentially M-T unliftable variety $B_i$
corresponding to this factor, we then have $r=1/2$.
This means that the construction of the intermediate group $H'$ can be
shunted in this case. 
It follows that the Mumford--Tate group of $B$ is contained in the
image of $\wtH^\der\times L^\times$ in $\GL(V\otimes_{K_0}W)$.

This argument is also valid for factors of type $A_k$ for
which all the numbers $r$ are equal to $1/2$. 
\subsubsection{}
Instead of using \cite{deligne:K3} to prove the fact $\tih$ defines a
polarizable Hodge structure on $V\otimes_{K_0}W$, one
may also explicitly construct a polarization. 
This is not difficult, using the autoduality of $V$ as 
representation of $\wtH^\der$. 
\subsection{Examples.}
\subsubsection{}
Let $A/\CC$ be an abelian variety arising from Mumford's construction,
see \cite{mumford:shimura}. 
In this case one has $G_{A/\bQ}^\der\cong\SL_2^3/\wtN$, where
\[
\wtN=\{(\epsilon_1,\epsilon_2,\epsilon_3)\mid 
\epsilon_i=\pm1\textrm{ for }i=1,2,3\textrm{ and } 
\epsilon_1\epsilon_2\epsilon_3=1\}. 
\]
The Mumford--Tate group of the M-T unliftable and M-T decomposed weak
M-T lift $B$ of $A$ satisfies $G_{B/\bQ}^\der\cong\SL_2^3$.
This is the example studied in detail in \cite{noot:lifting}.
\subsubsection{}
There exist simple abelian varieties $A/\CC$ for which $G_A^\ad$ is
not simple. 
One can construct such an example where 
$G_A^\der\cong G_1\times G_1/N$ with 
\[
G_{1/\bQ}\cong G_{2/\bQ}\cong\SL_2^2
\] 
and 
$N_{/\bQ}=\{(\pm1,\pm1)\}\subset\SL_2^2\cong G_i$ embedded diagonally into
$G_1\times G_2$. 
In this case, the M-T unliftable and M-T decomposed weak
M-T lift $B$ of $A$ is a product $B=B_1\times B_2$ and
$G_{B_i}^\der=G_i$. 
\subsubsection{}
There exist simple abelian varieties $A/\CC$ for which $G_A^\ad$ is
absolutely simple of type $D_k^\RR$, with $k$ even, $\cG_{\QQ}$ acting
trivially on the Dynkin diagram and where $V=\rH_\B^1(A(\CC),\QQ)$
decomposes over $\bQ$ as a multiple of the representation of
$G_A^\der$ with highest weight $\varpi_k$. 
Since $G_A^\der$ acts faithfully on $V$, it is not simply connected. 
Let $B$ be the M-T unliftable and M-T decomposed weak
M-T lift of $A$. 
Then $G_B^\der$ is the universal cover of $G_A^\der$ and $B$ is a
product $B\cong A\times B'$, where $W=\rH_\B^1(B'(\CC),\QQ)$
decomposes over $\bQ$ as a multiple of the representation of 
$G_B^\der$ with highest weight $\varpi_{k-1}$. 
\begin{cor}
For every abelian motive $M/\CC$ there exists an essentially
Mumford--Tate unliftable and Mumford--Tate decomposed abelian variety
$B/\CC$ which provides a weak Mumford--Tate lift for $M$. 
\end{cor}
\begin{proof}
There exist an abelian variety $A/\CC$ and a surjection of the
corresponding 
Mumford--Tate groups $G_A\rightarrow G_M$ commuting with the maps $h_A$
and $h_M$.
Let $B'/\CC$ be the essentially Mumford--Tate unliftable and
Mumford--Tate decomposed weak Mumford--Tate lift for $A$ provided by
theorem~\ref{unliftable_lift_exists}. 
This gives a
morphism $\pi_{B'}\colon G_{B'}^\der\rightarrow G_M^\der$. 
The fact that $B'$ is Mumford--Tate decomposed implies that 
there is an isogeny $B'\sim\prod_{i\in I}B_i$, where the $B_i$ are
M-T decomposed abelian varieties such that the groups $G_{B_i}^\der$
are $\QQ$-simple,
and such that $G_{B'}^\der\cong\prod_{i\in I} G_{B_i}^\der$. 
Let $J\subset I$ be the subset of indices $i$ such that
$G_{B_i}^\der$ is not in the kernel of $\pi_{B'}$ and let
$B=\prod_{i\in J} B_i$. 
Then $G_B=\prod_{i\in J} G_{B_i}$ and 
$\bigl(\pi_{B'}\bigr)_{|G^\der_B}$ is an
isogeny from $G_B^\der$ onto $G_M^\der$, so $B$ verifies the condition
of the corollary. 
\end{proof}
\section{Mumford--Tate liftings and motives}
\label{sect_descent_MT-adjoint}
\begin{prop}\label{descent_MT-adjoint}
Suppose that $A_{/\CC}$ and $B_{/\CC}$ are abelian varieties over $\CC$,
let $(G_A,h_A)$ and $(G_B,h_B)$ be the associated Mumford--Tate data
and assume that there exists an isomorphism $G_A^\ad\cong G_B^\ad$
such that $\pi_A^\ad\circ h_A=\pi_B^\ad\circ h_B$. 

Let $F\subset\CC$ be an algebraically closed field. 
Then there exists an abelian variety $A/F$ such
$A\otimes_F\CC=A_{/\CC}$ 
if and only if there exists an abelian variety $B/F$ such that
$B\otimes_F\CC=B_{/\CC}$.
\end{prop}
\begin{proof}
This is proved as in \cite[4.5]{noot:lifting}.
Let $h_A\colon S\rightarrow G_{A/\RR}$ and 
$h_B\colon S\rightarrow G_{B/\RR}$ be the maps defining the
Hodge structures on $\rH_\B^1(A(\CC),\QQ)$ and $\rH_\B^1(B(\CC),\QQ)$ 
respectively. 
Let $X_A$ and $X_B$ be the $G_A(\RR)$- and $G_B(\RR)$-conjugacy
classes of $h_A$ and $h_B$. 
The main theorem of \cite{deligne:shimura2} implies that 
for all compact open subgroups $K_A\subset G_A(\AA_f)$ and 
$K_B\subset G_B(\AA_f)$, one can construct quasi-canonical models 
${}_{K_A}M_{/\bQ}(G_A,X_A)^0$ and ${}_{K_B}M_{/\bQ}(G_B,X_B)^0$ 
over $\bQ$ of the corresponding connected Shimura varieties. 
Here $\bQ$ is the algebraic closure of $\QQ$ in $\CC$.

For $K_A$ and $K_B$ sufficiently small, there exist ``universal''
abelian schemes 
$\mathcal{A}\rightarrow{}_{K_A}M_{/\bQ}(G_A,X_A)^0$ and
$\mathcal{B}\rightarrow{}_{K_B}M_{/\bQ}(G_B,X_B)^0$ and
points 
\[
a\in{}_{K_A}M_{/\bQ}(G_A,X_A)^0(\CC),\;
b\in{}_{K_B}M_{/\bQ}(G_B,X_B)^0
\]
such that  $A=\mathcal{A}_a$ and $b=\mathcal{B}_b$. 

Let $G^\ad=G_A^\ad\cong G_B^\ad$. 
One can choose $K_A$, $K_B$ and $K^\ad\subset G^\ad(\AA_f)$
such that, in addition to the above conditions, there is a diagram
\[
\xymatrix{
{}_{K_A}M_{/\bQ}(G_A,X_A)^0 \ar[rd] &  & {}_{K_B}M_{/\bQ}(G_B,X_B)^0 \ar[ld] \\
& {}_{K^\ad}M_{/\bQ}(G^\ad,X^\ad)^0,
         }
\]
in which each arrow is a quotient map for the action of a finite group
and hence is a finite morphism 
(cf.~\cite[2.7.11(b)]{deligne:shimura2}) 
and such that $a$ and $b$ have
the same image in ${}_{K^\ad}M_{/\bQ}(G^\ad,X^\ad)(\CC)^0$. 
The proposition follows. 
\end{proof}
\begin{cor}\label{descent_strong_MT-lift}
The statement of the proposition is true in particular if $B_{/\CC}$ provides a
(weak) Mumford--Tate lift of $A_{/\CC}$. 
\end{cor}
\begin{prop}\label{MT-lift_motive}
Let $F\subset\CC$ be an algebraically closed field and
let $A$ and $B$ 
abelian varieties over $F$ such that $B_{/\CC}$ provides a Mumford--Tate
lift of $A_{/\CC}$. 
Then the motive $h^1(A)$ belongs to the category
$\langle h^1(B),\QQ(1)\rangle$.

The map $G_B\rightarrow G_A$ induced by the Betti realization of this
inclusion is the map given by the structure of $B$ as a Mumford--Tate
lift of $A$. 
\end{prop}
\begin{proof}
This generalises \cite[4.8--4.11]{noot:lifting}.

As explained in~\ref{prelim_betti}, 
it follows from \cite[6.25]{delmil:tannakian} that the Betti
realization induces an equivalence of 
$\langle h^1(B),\QQ(1)\rangle$ 
with 
$\langle\rH_\B^1(B(\CC),\QQ),\QQ(1)\rangle$, the tannakian subcategory 
of the category of Hodge structures generated by $\rH_\B^1(B(\CC),\QQ)$
and the Tate Hodge structure $\QQ(1)$. 
To prove the proposition, it therefore suffices to show that
the Hodge structure
$\rH_\B^1(A(\CC),\QQ)$ belongs to 
$\langle\rH_\B^1(B(\CC),\QQ),\QQ(1)\rangle$.
For the rest of the proof, we write $V_A=\rH_\B^1(A(\CC),\QQ)$ and
$V_B=\rH_\B^1(B(\CC),\QQ)$. 

By definition of the Mumford--Tate group, the underlying vector space
of any object of the category $\langle V_B,\QQ(1)\rangle$ of Hodge
structures naturally carries the
structure of a representation of $G_B$. 
This gives a $\otimes$-equivalence of $\langle V_B,\QQ(1)\rangle$
with the subcategory 
$\langle V_B,\QQ(1)\rangle_{\Rep}$ of $\Rep_\QQ(G_B)$.
As we saw in~\ref{prelim_betti}, for any object of $W$ the latter
category, the Hodge structure is given by composing the
morphism $h_B\colon S\rightarrow G_{B/\RR}$ with the action of $G_B$
on $W$. 

Let $(G_A,h_A)$ be the Mumford--Tate datum associated to $A$.
By hypothesis, there exists a central morphism 
$\pi\colon G_B\rightarrow G_A$ such that $h_A=\pi_\RR\circ h_B$. 
This makes every $\QQ$-linear representation $W$ of $G_A$ into a
representation of $G_B$ and, for any such $W$, it carries
the Hodge structure on $W$ defined by $h_A$ into the Hodge structure
defined by $h_B$. 
To prove the proposition it is therefore sufficient to show that $V_A$,
considered as representation of $G_B$, belongs to 
$\langle V_B,\QQ(1)\rangle_{\Rep}$.
This follows immediately, because $V_B$ is a faithful
representation of $G_B$ and $V_B\otimes\QQ(1)$ is
its dual, so one has 
\[
\langle V_B,\QQ(1)\rangle_{\Rep}=\Rep_\QQ(G_B),
\]
cf.~\cite[2.20]{delmil:tannakian} and its proof. 
\end{proof}
\begin{prop}\label{weak_MT-lift_motive}
Let $A$ and $B$ be abelian varieties over an algebraically closed
field $F\subset\CC$ such that $B_{/\CC}$ provides
a weak Mumford--Tate lift of $A_{/\CC}$. 
Then $h^1(A)$ belongs to $\langle h^1(B),\CM_F\rangle$. 

Taking the Betti realization, this inclusion induces a map between the
corresponding Mumford--Tate groups. 
On the derived group, this map is the map 
$\pi^\der\colon G_B^\der\rightarrow G_A^\der$ given by the structure of
$B$ as weak Mumford--Tate lift of $A$. 
\end{prop}
\begin{proof}
Let $T_A$ and $T_B$ be the connected components of the centres of
$G_A$ and of $G_B$ respectively. 
We write $V_A=\rH_\B^1(A(\CC),\QQ)$ and
$V_B=\rH_\B^1(B(\CC),\QQ)$ as in the proof of~\ref{MT-lift_motive} and
consider all spaces we encounter as representations of
$H=T_A\times T_B\times G_B^\der$. 
In the case of $V_A$, the group $H$ acts via 
\[
\begin{CD}
H @>>> T_A\times G_B^\der @>\id\times\pi^\der>> T_A\times G_A^\der
@>>> G_A
\end{CD}
\]
and in the case of $V_B$ it acts via 
$H\rightarrow T_B\times G_B^\der\rightarrow G_B$.

The groups $T_A$ and $T_B$ act on $V_A$ and $V_B$. 
Composing these representations with the projections 
$H\rightarrow T_A$ and $H\rightarrow T_B$ respectively we obtain
representations $V_A^c$ and $V_B^c$ of $H$. 
As $W=V_B\oplus V_A^c\oplus V_B^c$ is a faithful 
representation of $H$, it 
follows that $V_A$ belongs to the subcategory of
$\Rep_\QQ(H)$ $\otimes$-generated by $W$ and its dual. 
As there is an isomorphism of Hodge structures
$V_B^\vee\cong V_B\otimes\QQ(1)$, 
this in turn implies that $V_A$ belongs to the
tannakian subcategory of $\Rep_\QQ(H)$ generated by $V_B$ and 
the abelian representations of $H$, i.\ e.\ the representations
where $H$  acts through a commutative quotient. 

We fix an irreducible direct factor $W_2$ of $V_A$ in
$\Rep_\QQ(H)$.
There is an irreducible abelian $\QQ$-linear representation $W_3$
of $H$ such that $W_2$ is isomorphic to a subobject of 
$W_3\otimes V_B^{\otimes d}$.
Replace $V_B^{\otimes d}$ by an irreducible direct factor $W_1$ such
that $W_2$ still is a subobject of $W_3\otimes W_1$. 
The projection of $W_1\otimes W_1^\vee$ onto the trivial representation
induces a surjection 
$W_3\otimes W_1\otimes W_1^\vee\twoheadrightarrow W_3$. 
The composite of this surjection with the map 
\[
W_2\otimes W_1^\vee\hookrightarrow W_3\otimes W_1\otimes W_1^\vee
\]
deduced from the inclusion $W_2\hookrightarrow W_3\otimes W_1$ is a
map $W_2\otimes W_1^\vee\rightarrow W_3$. 
It is not difficult to check that this map is non-trivial and as 
$W_3$ was assumed to be irreducible, this implies that the map is
surjective.  
This proves that $W_3$ belongs to the tannakian subcategory of
$\Rep_\QQ(H)$ which is $\otimes$-generated by $W_1^\vee$ and $W_2$. 

It follows that the action of $H$ on $W_3$ factors through 
$G_A\times G_B$ and as $W_3$ is an abelian representation of $H$, it
is also an abelian representation of $G_A\times G_B$.
This means that $W_3$ carries a Hodge structure of CM-type and hence
that $V_A$ belongs to $\langle V_B,\HodgeCM\rangle$ as required. 
\end{proof}
\begin{cor}
Let $A$ and $B$ be abelian varieties over $F$ and let $G_A$ and $G_B$
be their respective Mumford--Tate groups. 
Assume that $B$ provides a weak Mumford--Tate lift of $A$ and that
$G_A^\der\cong G_B^\der$.
Then the categories 
$\langle h^1(A),\CM_F\rangle$ and $\langle h^1(B),\CM_F\rangle$
coincide. 
\end{cor}
\begin{cor}\label{AM_to_MTud-AV}
Let $F\subset\CC$ be an algebraically closed field and let
$M$ be an object of $\AV_F$. 
Then there exist an essentially Mumford--Tate unliftable and
Mumford--Tate decomposed abelian variety
$A$ over $F$ such that $M$ belongs to 
\[
\langle h^1(A),\CM_F\rangle.
\]
\end{cor}
\section{Motivic liftings of Galois representations}
\label{motivic_galois_liftings} 
\subsection{}\label{intro_galois-lifts}
Assume that $F\subset\CC$ is a number field, $M$ an abelian motive
over $F$ and $(G_M,h_M)$ the Mumford--Tate datum associated to $M$. 
As we recalled in~\ref{prelim_p-adic}, the 
$p$-adic Galois representation associated to the \'etale realization
of $M$ factors through a morphism 
\[
\rho_{M,p}\colon\cG_{F}\longrightarrow G_M(\Qp).
\]
There exists a finite extension $F'\supset F$ such that the
Mumford--Tate group of $M_{F'}$ is connected. 
The Mumford--Tate group of $M_{F'}$ is then equal to the connected
component of $G_M$. 
In what follows, we will assume that $G_M$ is already connected 
(replacing $F$ by a finite extension if necessary).

It follows from theorems of Tsuji, \cite{tsuji:cst} and De Jong,
\cite{dejong:alterations} that, for every
$p$, the representation $\rho_{M,p}$ of $\cG_F$ on $\rH_\et(M,\Qp)$ is
geometric in the sense of 
Fontaine and Mazur, \cite[\S1]{fonmaz}. 
Here a representation of $\cG_F$ on a finite dimensional $\Qp$-vector
space is called \emph{geometric} if 
\begin{itemize}
\item
it is unramified outside a finite set of non-archimedean places of $F$
and 
\item
for each valuation $\bar v$ of $\bar F$ of residue
characteristic $p$, the restriction to the inertia
group $\cI_{F,\bar v}$ is potentially semi-stable
(cf.~\cite{fontaine:repr_semist}). 
\end{itemize}
More generally, for a linear algebraic group $G$ over $\Qp$ and a continuous
morphism $\rho_p\colon\cG_{F}\longrightarrow G(\Qp)$, we will
say that $\rho_p$ is \emph{geometric} if there exists a faithful
representation $V_p$ of $G$ such that the resulting
representation of $\cG_F$ on $V_p$ is geometric. 
This is the case if and only if the representation of $\cG_F$ on $W_p$
is geometric for any representation $W_p$ of $G$.
\begin{dfn}\label{def_geometric_lift}
Assume that $G$ and $\wtG$ are linear algebraic groups over $\QQ$ or
over $\Qp$ and that
$\rho_p\colon\cG_F\rightarrow G(\Qp)$ and
$\trho_p\colon\cG_F\rightarrow\wtG(\Qp)$ are geometric Galois
representations. 
We will say that $\trho_p$ is a \emph{geometric lift} of $\rho_p$ if
there exists a central isogeny $\pi\colon\wtG\rightarrow G$ such
that $\rho_p=\pi\circ\trho_p$. 
If $\rho_p$ does not admit any geometric lift with $\ker\pi$
non-trivial, it will be called \emph{geometrically unliftable}. 

If there exists a central isogeny 
$\pi^\der\colon\wtG^\der\rightarrow G^\der$ such that there is an equality
$\tilde\pi^\ad\circ\trho_p=\pi^\ad\circ\rho_p$ then $\trho_p$ is
said to be a \emph{weak geometric lift} of $\rho_p$. 
As in~\ref{def_MT-lift}, the maps $\pi^\ad$ and $\tilde\pi^\ad$ are
the projections $G\rightarrow G^\ad$ and $\wtG\rightarrow\wtG^\ad$. 
The central isogeny $\pi^\der$ induces an isomorphism 
$G^\ad\cong\wtG^\ad$, giving a sense to the equality 
$\tilde\pi^\ad\circ\trho_p=\pi^\ad\circ\rho_p$.
We will say that $\rho_p$ is \emph{essentially geometrically
unliftable} if it does not admit a weak geometric lift with
$\ker\pi^\der$ non-trivial. 
\end{dfn}
\begin{prop}\label{MT_lift-galois_lift}
Suppose that $F\subset\CC$ is a number field, $A/F$ an abelian
variety with connected Mumford--Tate group $G_A$ and that $B$
is an abelian variety over a finite extension $F'\supset F$ such that
$B_{/\CC}$ provides a Mumford--Tate lift of $A_{/\CC}$. 

The number field $F'$ can be chosen such that, for every prime number
$p$, the 
Galois representation $\rho_{A,p}$ of $\cG_{F'}$ on 
$\rH_\et^1(A_{\bar F},\Qp)$
belongs to the subcategory of $\Rep_{\Qp}(\cG_{F'})$ 
$\otimes$-generated by $\rH_\et^1(B_{\bar F},\Qp)$ and $\Qp(1)$.
For every $p$, the morphism $G_B\rightarrow G_A$ then realises 
the Galois representation 
$\rho_{B,p}$ as a geometric lift of $\rho_{A,p}$.
\end{prop}
\begin{proof}
By proposition~\ref{MT-lift_motive}, the motive $h^1(A_{\bar F})$
belongs to $\langle h^1(B_{\bar F},\QQ(1)\rangle$ and, for $F'$ large
enough, this is already the case over $F'$. 
Taking the
$p$-adic \'etale realizations, this implies the corresponding
statement for the Galois representations. 
The inclusion of $h^1(A_{F'})$ in 
$\langle h^1(B_{F'}),\QQ(1)\rangle$
corresponds to a morphism $G_B\rightarrow G_A$ and taking the $p$-adic
\'etale realizations this gives rise to a commutative diagram 
\[
\xymatrix{&& G_B(\Qp)\ar[d]\\
\cG_{F'}\ar[rru]^{\rho_{B,p}}\ar[rr]^{\rho_{A,p}}&&G_A(\Qp)
	}
\]
proving the proposition. 
\end{proof}
\subsection{}
As above, assume that $F$ is a number field contained in $\CC$. 
Following the notations of \cite[\S6]{fonmaz}, let
$\RepCM{\Qp}{(\cG_F)}$, or $\RepCM{}{}$ if no confusion is likely,
denote the tannakian subcategory of 
$\Rep_{\Qp}(\cG_F)$ consisting of the potentially abelian geometric
representations, in other words, the geometric representations such that the
restriction to a subgroup of $\cG_F$ of finite index has abelian
image. 
It follows from \cite[\S6]{fonmaz} that $\RepCM{}{}$ is the tannakian
subcategory of $\Rep_{\Qp}(\cG_F)$ generated by the representations 
factoring through 
finite groups and the representations of the form 
$\rH_\et^1(A_{\bar F},\Qp)$ for $A/F$ an abelian variety which is
potentially of CM-type. 
Thus, $\RepCM{}{}$ is the tannakian category of the $p$-adic \'etale
realizations of the objects of $\CM_F$. 
\begin{prop}\label{weak_MT_lift-galois_lift}
Suppose that $F\subset\CC$ is a number field, $A/F$ an abelian
variety with connected Mumford--Tate group $G_A$ and that $B$
is an abelian variety over a finite extension $F'\supset F$ such that
the Mumford--Tate group 
$G_B$ is connected and $B_{/\CC}$ provides a 
weak Mumford--Tate lift of $A_{/\CC}$. 

For every prime number
$p$, the representation $\rho_{A,p}$ of $\cG_{F'}$ on
$\rH_\et^1(A_{\bar F},\Qp)$ is an object of the subcategory 
$\langle\rH_\et^1(B_{\bar F},\Qp),\RepCM{\Qp}{(\cG_{F'})}\rangle$
of $\Rep_{\Qp}(\cG_{F'})$.
Via the map $G_B^\der\rightarrow G_A^\der$, the Galois
representation $\rho_{B,p}$ provides a
weak geometric lift of $\left(\rho_{A,p}\right)_{|\cG_{F'}}$.
\end{prop}
\begin{proof}
It follows from proposition~\ref{weak_MT-lift_motive} that there is a
finite extension $F''$ of $F'$ such that the motive $h^1(A_{F''})$ 
belongs to $\langle h^1(B_{F''}),\CM_{F''}\rangle$. 
As the Mumford--Tate group of $\langle h^1(B_{F''}),\CM_{F''}\rangle$ is
a pro-algebraic group with derived group $G_B^\der$, this gives rise
to an isogeny $\pi\colon G_B^\der\rightarrow G_A^\der$ and $\pi$
induces an isomorphism $\pi^\ad\colon G_B^\ad\rightarrow G_A^\ad$.
Taking $p$-adic \'etale realizations, there is a commutative diagram
\begin{equation}\label{wgl-diagram}
\parbox[c]{5cm}{\xymatrix{
 &&G_B(\Qp)\ar[dr]\\
 \cG_{F''}\ar[rru]^{\rho_{B,p}}\ar[rrd]^{\rho_{A,p}}&&&
  G_A^\ad(\Qp)=G_B^\ad(\Qp).\\
  &&G_A(\Qp)\ar[ur]
	 }}
\end{equation}
This proves all statements of the proposition with $F''$ instead of
$F'$. 

Faltings' theorem, \cite[Satz~4]{faltings:endlichkeit} implies that
if $C$ is the commuting algebra of $G_A$ in
$\End(\rH_B^1(A(\CC),\QQ))$, then $C\otimes\Qp$ is the commuting algebra of 
$\rho_{A,p}(\cG_{F''})$ in $\End(\rH_\et^1(A(\CC),\Qp))$. 
This implies that that the centralizer of the image of $\cG_{F''}$ in
$G_A^\ad(\Qp)$ is trivial. 
It follows from lemma~\ref{wgl_after_ext} that the
diagram~\ref{wgl-diagram} also commutes with $\cG_{F''}$ replaced by 
$\cG_{F'}$.
Lemma~\ref{if_wgl_then_tannakian} implies that the representation of
$\cG_{F'}$ on 
$\rH_\et^1(A_{\bar F},\Qp)$ is an object of 
$\langle\rH_\et^1(B_{\bar F},\Qp),\RepCM{\Qp}{(\cG_{F'})}\rangle$. 
\end{proof}
\begin{lemma}\label{wgl_after_ext}
Assume that $p$ is a prime number and that $G$ is a connected linear algebraic
group over $\Qp$. 
Let $F$ be a number field and  let 
$\rho_1,\rho_2\colon\cG_F\rightarrow G(\Qp)$ be Galois representations
coinciding on $\cG_{F'}$ for some finite extension
$F'$ of $F$. 
Also assume that the centralizer of 
$\rho_1(\cG_{F'})=\rho_2(\cG_{F'})$ in $G$ is
trivial. 
Then $\rho_1=\rho_2$ on $\cG_F$.
\end{lemma}
\begin{proof}
It is sufficient to treat the case where $F'$ is a Galois extension of
$F$. 

Let $\delta\colon\cG_F\rightarrow G(\Qp)$ be defined by
$\delta(\sigma)=\rho_1(\sigma)\rho_2(\sigma)^{-1}$. 
This map satisfies the cocycle condition 
$\delta(\sigma\tau)=
\delta(\sigma)\left(\rho_2(\sigma)\delta(\tau)\rho_2(\sigma)^{-1}\right)$.
This implies that $\delta$ is constant on the classes $\sigma\cG_{F'}$
for $\sigma\in\cG_F$. 
As $\cG_{F'}$ is normal, it follows that $\delta$ is also constant
on the classes $\cG_{F'}\sigma$ and it follows that for all
$\sigma\in\cG_F$ and $\tau\in\cG_{F'}$ one has 
$\delta(\sigma)=\rho_2(\tau)\delta(\sigma)\rho_2(\tau)^{-1}$. 
Therefore 
$\delta(\sigma)$ lies in the centralizer of $\rho_2(\cG_{F'})$ and we
conclude that $\delta$ is trivial. 
\end{proof}
\begin{lemma}\label{if_wgl_then_tannakian}
Let $p$ be a prime number, $G_1$ and $G_2$ connected linear algebraic
groups over $\Qp$ and $\pi^\der\colon G_1^\der\rightarrow G_2^\der$ a
central isogeny. 
Let $V_1$ be a faithful $\Qp$-linear representation of
$G_1$ and let $V_2$ be any $\Qp$-linear representation of $G_2$. 

Let $F$ be a number field and  let 
$\rho_i\colon\cG_F\rightarrow G_i(\Qp)$ (for $i=1,2$) be geometric Galois
representations.
Assume  that $V_1^\vee$ lies in the subcategory of
$\Rep_{\Qp}(\cG_{F})$  $\otimes$-generated by $V_1$ and
$\RepCM{\Qp}{(\cG_{F})}$ 
and that $\pi_1^\ad\circ\rho_1=\pi_2^\ad\circ\rho_2$,
where the $\pi_i^\ad\colon G_i\rightarrow G_i^\ad$ are the canonical
projections and $G_1^\ad=G_2^\ad$ is the identification induced
by $\pi^\der$. 
Then $V_2$ is an object of the subcategory of
$\Rep_{\Qp}(\cG_{F})$  $\otimes$-generated by $V_1$ and
$\RepCM{\Qp}{(\cG_{F})}$. 
\end{lemma}
\begin{proof}
Fix an object $V_3$ of $\RepCM{\Qp}(\cG_{F})$ such that $V_1^\vee$ lies
in the subcategory of $\Rep_{\Qp}(\cG_{F})$  $\otimes$-generated by
$V_1$ and $V_3$. 
Let $T_3$ be the Zariski closure of the image of $\cG_F$ in $\GL(V_3)$
and let $\rho_3\colon\cG_F\rightarrow T_3(\Qp)$ be the morphism giving
the action of $\cG_F$ on $V_3$. 

Let $T_1$ and $T_2$ be the connected components of the centres of
$G_1$ and of $G_2$ respectively. 
Put $H=T_1\times T_2\times T_3\times G_1^\der$ and consider the natural maps 
$H\rightarrow G_1$ and $H\rightarrow G_2$. 
Via these maps, we consider $V_1$ and $V_2$ as representations
of $H$. 
It can be shown exactly as in the proof of
proposition~\ref{weak_MT-lift_motive} that $V_2$ belongs to the
subcategory of $\Rep_{\Qp}(T_3\times G_A\times G_B)$ generated by
$V_1$ and the abelian representations. 
For any abelian representation $W_3$ of $T_3\times G_A\times G_B$, the induced
Galois representation belongs to $\RepCM{}{(\cG_F)}$ so this implies
that $V_2$ belongs to the subcategory 
$\langle V_1,\RepCM{}{(\cG_F)}\rangle$ of $\Rep_{\Qp}(\cG_{F})$.
\end{proof}
\subsection{}\label{intro_non_spinorial_lift}
We keep the above notations, i.\ e.\ $F\subset\CC$ is a number field,
$A/F$ an abelian variety and $G_A$ its Mumford--Tate group, which is
assumed to be connected.
For each
prime number $p$, we denote by $\rho_{A,p}\colon\cG_F\rightarrow G_A(\Qp)$ the
$p$-adic Galois representation associated to $A$.

Fix an algebraic group $\wtG$ over $\QQ$ and a central isogeny 
$\wtG^\der\rightarrow G_A^\der$ and suppose that for every $p$ in a
set $P$ of prime numbers $\trho_p\colon\cG_F\rightarrow\wtG(\Qp)$ is a
weak geometric lift  of $\rho_{A,p}$. 
\begin{thm}\label{thm_non_spinorial_lift_from_AM}
Assume that any simple factor of $\wtG^\der_{/\CC}$ lying over a
factor of $G_{A/\CC}^\der$ of type $D_k^\HH$ is a quotient of the
$h$-maximal cover, cf.~\ref{D_k-factors}. 
Let $F'$ be a finite extension of $F$ and $B$
an essentially Mumford--Tate unliftable abelian variety over $F'$ with
connected Mumford--Tate group $G_B$ such that $B_{/\CC}$
provides a weak Mumford--Tate lift of $A_{/\CC}$. 
For each prime number $p$, let $V_{B,p}=\rH^1_\et(B_{\bar F},\Qp)$
be the $p$-adic representation of $\cG_{F'}$ associated to $B$. 

Then the map $G_B^\der\rightarrow G_A^\der$ lifts to
$G_B^\der\rightarrow\wtG^\der$. 
For every $p\in P$ 
\begin{itemize}
\item
the
representation $\rho_{B,p}$ is a weak 
geometric lift of the restriction $\left(\trho_p\right)_{|\cG_{F'}}$
and 
\item
for every
representation $\wtV_p$ of $\wtG_{/\Qp}$, the Galois representation on
$\wtV_p$ is an object of 
$\langle V_{B,p},\RepCM{}{}\rangle$.
\end{itemize}
\end{thm}
\begin{proof}
We fix a faithful self-dual representation $\wtV$ of $\wtG$ and for
each prime number $p$ we write $\wtV_p=\wtV\otimes_\QQ\Qp$. 
As $\wtV$ generates the tannakian category of representations of
$\wtG$, it is sufficient to prove the corollary for the
representations $\wtV_p$. 

Since $B$ is essentially Mumford--Tate unliftable, it follows from
theorem~\ref{MT-group_unliftable-AV} that $G_B^\der$ is $h$-maximal. 
It follows that the map $G_B^\der\rightarrow G_A^\der$ lifts to a map
$G_B^\der\rightarrow\wtG^\der$. 
Let $\rho_{B,p}\colon\cG_{F'}\rightarrow G_B(\Qp)$ be the map giving the
Galois representation on $V_{B,p}$.
Write $\pi_A^\ad\colon G_A\rightarrow G_A^\ad$, 
$\pi_B^\ad\colon G_B\rightarrow G_B^\ad$ and 
$\tilde\pi^\ad\colon\wtG\rightarrow\wtG^\ad$ for the projections. 
For any $p\in P$, the
proposition~\ref{weak_MT_lift-galois_lift} and the fact that $\trho_p$
is a weak geometric lift of $\rho_{A,p}$ imply that 
$\tilde\pi^\ad\circ\trho_{p|\cG_{F'}}=
\pi_A^\ad\circ\rho_{A,p|\cG_{F'}}=
\pi_B^\ad\circ\rho_{B,p}$. 
The remaining statement of the theorem follows from
lemma~\ref{if_wgl_then_tannakian}. 
\end{proof}
\subsection{Important remark.}\label{B'_in_thm_exists}
Concerning the abelian variety $B/F'$ which appears in
propositions~\ref{MT_lift-galois_lift}
and~\ref{weak_MT_lift-galois_lift} and in
theorem~\ref{thm_non_spinorial_lift_from_AM}, 
it follows from theorem~\ref{unliftable_lift_exists} and
proposition~\ref{descent_MT-adjoint} that there exist a number field
$F'$ and an essentially M-T unliftable and M-T decomposed weak
Mumford--Tate lift 
$B/F'$ of $A$ as in the propositions and in the theorem.
The condition that $G_B$ is connected can be forced by replacing $F'$
by a finite extension. 
\begin{cor}\label{cor_non_spinorial_lift_from_AM}
Let notations be as in~\ref{intro_non_spinorial_lift} with
$\wtG^\der\rightarrow G_A^\der$ satisfying the hypotheses of the
theorem. 
Then, for every $p\in P$ and every representation $\wtV_p$ of
$\wtG_{/\Qp}$, the induced representation of $\cG_F$ on $\wtV_p$
occurs in the $p$-adic \'etale realization of an abelian motive. 
\end{cor}
\begin{proof}
The remark and the theorem imply that there is is a finite extension
$F'$ of $F$ 
such that the representation of $\cG_{F'}$ on $\wtV_p$ occurs in the
$p$-adic \'etale realization of an object $M'$ of $\AV_{F'}$.
The representation of $\cG_F$ on $\wtV_p$ then occurs in the $p$-adic
\'etale realization of the Weil restriction $\Res_{F'/F}M'$, which is
also in $\AV_F$. 
\end{proof}
\section{Abelian varieties with Mumford--Tate group of type $D_k^\HH$}
\label{AV-D_k^H-section}
\subsection{}\label{DkH_intro}
The following notation and hypotheses will be in force until
definition~\ref{def_lifted_abelian_DkH}. 
We let $A/\CC$ be a simple abelian variety and $(G_A,h_A)$ the associated
Mum\-ford--Tate datum. 
This implies that $G_A$ is connected. 
We assume throughout that $A$ is essentially Mumford--Tate unliftable,
Mumford--Tate decomposed and that $G_A$ is of type $D_k^\HH$ with
$k\geq4$.

It follows that there exist a totally
real number field $K_0$ and an absolutely simple algebraic group
$G^s/K_0$ such that $G_A^\der=\Res_{K_0/\QQ}G^s$.
By assumption, the group $G^s/K_0$ is $h$-maximal in the sense of
\ref{D_k-factors}. 
The representation of $G_A^\der$ on $\rH_\B^1(A(\CC),\QQ)$
decomposes over $\bQ$ as a multiple of the direct sum of the standard
(orthogonal) 
representations of the different factors of $G_{A/\bQ}^\der$.
As in \cite{deligne:shimura2} and in \S\ref{sect_MT_liftings} of this
paper, denote 
by  $I=\{\iota\colon K_0\hookrightarrow\CC\}$ the set of complex
embeddings of $K_0$. 
As $K_0$ is totally real, $I$ is also the set of real embeddings of
$K_0$.
The Dynkin diagram of $G_{A/\CC}^\der$ is a disjoint union, indexed by
$I$, of diagrams of type $D_k$.
The Hodge cocharacter $\mu_A\colon\Gm{/\CC}\rightarrow G_{/\CC}$
associated to $A$ projects trivially on some factors of $G^\ad_{/\CC}$. 
On the other factors, the conjugacy class of the projection 
is dual to one of
the vertices $\alpha_{k-1}$ or $\alpha_k$ (or possibly $\alpha_1$ if
$k=4$) of the corresponding
component of the Dynkin diagram. 
Without loss of generality, we will henceforth assume that, on these
factors, it corresponds to $\alpha_k$.
As in section~\ref{sect_MT_liftings},
let $I_c\subset I$ be the set of embeddings corresponding to the
factors onto which $\mu_A$ projects trivially and let $I_{nc}=I-I_c$.
Recall that $I_c$ is the set of embeddings 
$\iota\colon K_0\hookrightarrow\RR$ such that the factor of
$G_{A/\RR}^\der$ corresponding to $\iota$ is compact and $I_{nc}$ is
the set of real embeddings of $K_0$ for which the corresponding factor of
$G_{A/\RR}^\der$ is non-compact.

As the Hodge filtration and its complex conjugate are opposite
filtrations, the complex conjugate of $\mu$ is conjugate to
$\mu^{-1}$, up to a central cocharacter. 
This implies that complex conjugation acts on the Dynkin diagram by
the main involution and hence that it acts 
trivially 
if $k$ is even and exchanges $\alpha_{k-1}$ and $\alpha_k$ on
every factor if $k$ is odd. 
It follows that $\cG_\QQ$ (or $\Aut(\CC)$) acts on the Dynkin
diagram through 
$\Gal(K/\QQ)$ for a number field $K\supset K_0$ with $[K:K_0]=1$ or
$2$ and which is totally real if $k$ is even, a CM field if $k$ is
odd. 
In particular $[K:K_0]=2$ if $k$ is odd. 
Note that the statement is also true for $k=4$, because it follows from the
definition of the case $D_4^\HH$ (see~\ref{D_k-factors}) that the
vertices $\alpha_1$ of the connected components of the Dynkin diagram
form a $\cG_\QQ$-orbit, so the stabilizer of a connected component is
of order at most $2$. 
\subsection{Construction of $\wtG$ and $\tmu$.}\label{constr_tildes}
We aim to construct an algebraic group $\wtG/\QQ$ such that
$\wtG^\der$ is simply connected and that $\wtG^\ad=G_A^\ad$, together
with a cocharacter $\tmu$ of $\wtG_{/\CC}$ such that
$\tilde\pi^\ad\circ\tmu=\pi^\ad\circ\mu$. 
The argument is strongly inspired by \cite{deligne:shimura2}, see also
the proof of theorems~\ref{MT-group_unliftable-AV}
and~\ref{unliftable_lift_exists}. 
Here, as before, $\pi^\ad$ and $\tilde\pi^\ad$ are the projections
$G_A\rightarrow G_A^\ad$ and $\wtG\rightarrow G_A^\ad$ respectively. 

Let $\wtG^s/K_0$ be the simply connected cover of $G^s$. 
Consider the direct sum of the representations of $\wtG^s_{/\bQ}$
with highest weights $\varpi_{k-1}$ and $\varpi_k$.
A multiple of that representation can be defined over $K_0$, 
let $W^s$ be the resulting representation of $\wtG^s$.
By construction, there is a decomposition
$W^s\otimes_{K_0}K=W^s_1\oplus W^s_2$, where $W^s_1$ (resp.\
$W^s_2$) is a multiple of the representation with highest weight 
$\varpi_{k-1}$ (resp.\ $\varpi_k$). 
A non-trivial element of $\Gal(K/K_0)$ exchanges the factors $W^s_1$
and $W^s_2$ of this decomposition. 
If $[K:K_0]=2$, then the composite map 
$W^s\subset W^s\otimes_{K_0}K\twoheadrightarrow W^s_1$ is an
isomorphism of $K_0$-vector spaces and endows $W^s$ with a
structure of $K$-vector space. 

Put $\wtG^\der=\Res_{K_0/\QQ}\wtG^s$ and let $W$ be the
rational representation of $\wtG^\der$ deduced from $W^s$. 
Since $W$ is the underlying $\QQ$-vector space of $W^s$, it
carries a structure of $K$-vector space. 
The cocharacter $\pi^\ad\circ\mu$ of $G^\ad_{/\CC}$ lifts to a
quasi-cocharacter $\nu$ of $\wtG^\der_{/\CC}$, cf.~\ref{prel_basic}. 
\begin{lemma}\label{weights_nu}
There is a decomposition 
\begin{equation}\label{dec_WC-K0}
W\otimes_\QQ\CC\cong\bigoplus_{\iota\in I}
\left(W^{(\iota)}_1\oplus W^{(\iota)}_2\right) 
\end{equation}
such that, for each $\iota\in I$, the representation $W^{(\iota)}_1$
(resp. $W^{(\iota)}_2$)
is a multiple of the irreducible representation with highest weight
$\varpi_{k-1}$ (resp. $\varpi_k$) of the factor of $G_{A/\CC}^\der$
corresponding to $\iota$. 

The weights of $\nu$ on the $W^{(\iota)}_j$ are trivial for $j=1,2$ if 
$\iota\in I_c$. 
For $\iota\in I_{nc}$, the weights of $\nu$ on $W^{(\iota)}_1$ are of
the form $\frac{k-2}{4}-m_1$ and those on $W^{(\iota)}_2$ are of the
form $\frac{k}{4}-m_2$ with $m_1,m_2\in\ZZ$. 
If $k$ is even $m_1$ runs from $0$ to $\frac{k-2}{2}$ and $m_2$ from $0$
to $\frac{k}{2}$. 
For $k$ odd, both $m_1$ and $m_2$ run from $0$ to $\frac{k-1}{2}$.
\end{lemma}
\begin{proof}
The decomposition of $W^s$ induces the
decomposition~(\ref{dec_WC-K0}). 
The quasi-cocharacter $\nu$ projects
trivially to the factors of $\wtG^\der_{/\CC}$ corresponding to the
$\iota\in I_c$ and non-trivially to the other factors. 
The highest weights of $\nu$ on $W^{(\iota)}_1$ and $W^{(\iota)}_2$  
can easily be deduced from the information collected in
\cite[Table~1.3.9]{deligne:shimura2}, it is the rational number corresponding
to $\varpi_{k-1}$ resp.\ $\varpi_k$ in that table. 
The lowest weight of $\nu$ on $W^{(\iota)}_1$ 
is the opposite of the number corresponding to $\varpi_{k-1}$ if $k$ is
even and the opposite of the number corresponding to $\varpi_k$ if $k$ is 
odd.
For $W^{(\iota)}_2$, the converse is the case. 
\end{proof}
We continue the construction of $\wtG$ and $\tmu$. 
The case where $k$ is even and the case where $k$ is odd will be
treated separately. 
\subsection{The case where $k$ is even.}\label{case_k_even}
In this case, $K$ is totally real and either $K=K_0$ or $[K:K_0]=2$. 
Let $L'$ be a totally imaginary quadratic extension of $K_0$, put
$L=KL'$ and define 
the algebraic torus $T_L$ over $\QQ$ by 
$T_L=\ker(N_{L/K})\subset L^\times$, where 
$N_{L/K}\colon L^\times\rightarrow K^\times$ 
is the field norm. 
The group $T_L$ naturally acts on $L$ and this gives rise to a
$\QQ$-linear representation $V_L$ of $T_L$. 
We have a natural structure of
$K$-vector space on $V_L$. 
For each embedding $\iota\colon K\hookrightarrow\CC$ we choose a
complex embedding of $L$ extending $\iota$ and this gives an
identification 
\[
T_{L/\CC}=\bigoplus_{\iota\in\wtI}\Gm{/\CC},
\]
where $\wtI$ is the set of complex embeddings of $K$. 

From now on, we will further distinguish the cases where $[K:K_0]=1$
and where $[K:K_0]=2$. 
\par\noindent\textbf{The subcase where $K=K_0$.}
Consider the decomposition~(\ref{dec_WC-K0}) of $W\otimes_\QQ\CC$.
Since $W^s=W_1^s\oplus W_2^s$ is a decomposition of $W^s$ as a direct
sum of $K_0$-vector
spaces, we obtain a $\QQ$-linear decomposition $W=W_1\oplus W_2$.
For $\iota\in I_{nc}$, the weights of $\nu$ on
$W^{(\iota)}_1$ are in $\frac{1}{2}+\ZZ$ if
$k\equiv0\pmod4$ and in $\ZZ$ if $k\equiv2\pmod4$.
For the weights on $W^{(\iota)}_2$, the converse is the case. 

Let $T=T_L\times T_L$ and define $W_3$ (resp.\ $W_4$) be the
representation of $T$ given by the action of the first (resp.\ the
second) factor $T_L$ on $V_L$. 
We define a quasi-cocharacter 
\[
\nu_L\colon\Gm{/\CC}\rightarrow T_{/\CC}=
\bigoplus_{\iota\in I}\Gm{/\CC}^2,
\]
by
\[
\nu_L(z)_\iota=
\left\{\begin{array}{ll} 
(1,1)&\textrm{if }\iota\in I_c\\
(\sqrt{z},1)&\textrm{if }\iota\in I_{nc}\textrm{ and }k\equiv0\pmod4\\
(1,\sqrt{z})&\textrm{if }\iota\in I_{nc}\textrm{ and }k\equiv2\pmod4,
\end{array}\right.
\]
where $\nu_L(z)_\iota$ is the component of $\nu_L(z)$ in the factor of
$T_{/\CC}$ corresponding to $\iota\in I$. 
Finally, let $\wtV$ be the representation of $\wtG^\der\times T$
defined by 
\[
\wtV=W_1\otimes_{K_0}W_3\oplus W_2\otimes_{K_0}W_4
\]
and let $\wtG$ be the image of $\wtG^\der\times T$ in $\GL(\wtV)$. 
This will not cause any confusion, since $\wtG^\der$ acts
faithfully on $\wtV$, so the derived group
of $\wtG$ is $\wtG^\der$.
The weights of $(\nu,\nu_L)$ on $\wtV$ are in $\ZZ$, so the projection
of the quasi-cocharacter $(\nu,\nu_L)$ from 
$(\wtG^\der\times T)_{/\CC}$ to
$\wtG_{/\CC}$ is a true cocharacter $\tmu$ of $\wtG_{/\CC}$. 
By construction,
\[
\tilde\pi^\ad\circ\tmu=\tilde\pi^\ad\circ\nu=\pi^\ad\circ\mu,
\]
so we have constructed the couple $(\wtG,\tmu)$ in this particular
subcase. 
\par\noindent\textbf{The subcase where $[K:K_0]=2$.}
The structure of $K$-vector space on $W$ gives rise to a decomposition 
\begin{equation}\label{dec_WC-K}
W\otimes_\QQ\CC\cong\bigoplus_{\tilde\iota\colon K\hookrightarrow\CC}
W^{(\tilde\iota)}.
\end{equation}
Recall that $\wtI$ is the set of complex embeddings of $K$.
Put $\wtI_c=\{\tilde\iota\in\wtI\mid\tilde\iota_{|K_0}\in I_c\}$ and
$\wtI_{nc}=\{\tilde\iota\in\wtI\mid\tilde\iota_{|K_0}\in I_{nc}\}$.
For $\tilde\iota\in\wtI_c$, the weights of $\nu$ on
$W^{(\tilde\iota)}$ are all $0$, for $\tilde\iota\in\wtI_{nc}$, these
weights are either in $\ZZ$ or in $\frac{1}{2}+\ZZ$. 
Let $\wtI_{nc,1}$ be the set of $\tilde\iota\in\wtI_{nc}$ where the former
possibility occurs and $\wtI_{nc,2}$ its complement in $\wtI_{nc}$. 
Note that for any $\iota\in I_{nc}$ one of the embeddings
$K\hookrightarrow\CC$ restricting to $\iota$ lies in $\wtI_{nc,1}$ and
the other one lies in $\wtI_{nc,2}$. 

This time we define a quasi-cocharacter 
\[
\nu_L\colon\Gm{/\CC}\rightarrow T_{L/\CC}=
\bigoplus_{\tilde\iota\in\wtI}\Gm{/\CC},
\]
by
\[
\nu_L(z)_{\tilde\iota}=
\left\{\begin{array}{ll} 
1&\textrm{if }\tilde\iota\in\wtI_c\cup\wtI_{nc,1}\\
\sqrt{z}&\textrm{if }\tilde\iota\in\wtI_{nc,2}.
\end{array}\right.
\]
As above, $\nu_L(z)_{\tilde\iota}$ is the component of $\nu_L(z)$ in
the factor of $T_{L/\CC}$ corresponding to $\tilde\iota\in\wtI$. 
Finally, let $\wtV=W\otimes_K V_L$ as representation of $\wtG^\der\times T$
and let $\wtG$ be the image of $\wtG^\der\times T$ in $\GL(\wtV)$. 
Once again, there is no risk of confusion, because the derived group
of $\wtG$ is $\wtG^\der$.
Since the weights of $(\nu,\nu_L)$ on $\wtV$ are in $\ZZ$, the projection
of $(\nu,\nu_L)$ from $(\wtG^\der\times T)_{/\CC}$ to $\wtG_{/\CC}$ is
a true cocharacter $\tmu$ of $\wtG_{/\CC}$. 
By construction,
$\tilde\pi^\ad\circ\tmu=\tilde\pi^\ad\circ\nu=\pi^\ad\circ\mu$,
so this achieves the construction of the couple $(\wtG,\tmu)$ in the
case where $k$ is even.
\subsection{The case where $k$ is odd.}\label{case_k_odd}
In this case, $K$ is a totally imaginary quadratic extension of
$K_0$. 
In order to unify this case with the previous one as much as possible,
we put $L=K$.
Let $T_L$ be the $\QQ$-algebraic torus $\ker(N_{L/K_0})$, where 
$N_{L/K_0}\colon L^\times\rightarrow K_0^\times$ is the field norm. 
As in the beginning of~\ref{case_k_even}, the action of $T_L$ on $L$
by multiplication on the left gives rise to a representation $V_L$
of $T_L$. 

We again consider the decomposition~(\ref{dec_WC-K}) and the subsets
$\wtI_c$ and $\wtI_{nc}$ of the set $\wtI$ of complex embeddings of
$K=L$. 
By lemma~\ref{weights_nu}, the weights of $\nu$ on $W^{(\tilde\iota)}$
are either in $\frac{1}{4}+\ZZ$ or in $-\frac{1}{4}+\ZZ$ for
$\tilde\iota\in\wtI_{nc}$ and trivial for $\tilde\iota\in\wtI_c$. 
Let $\wtI_{nc,1}$ the set of $\tilde\iota\in\wtI_{nc}$ for which the
weights of $\nu$ on $W^{(\tilde\iota)}$ are in $\frac{1}{4}+\ZZ$
and let $\wtI_{nc,2}$ be the complement of $\wtI_{nc,1}$ in $\wtI_{nc}$.
For each $\iota\in I_{nc}$, one of the embeddings
$K\hookrightarrow\CC$ restricting to $\iota$ lies in $\wtI_{nc,1}$ and
the other one lies in $\wtI_{nc,2}$. 
Thus, the map $r\colon\wtI\rightarrow I$ given by
$\tilde\iota\mapsto\tilde\iota_{|K_0}$ induces a
bijection of $\wtI_{nc,1}$ with $I_{nc}$. 
For $I_c$, we arbitrarily fix a subset $\wtI_{c,1}\subset\wtI_c$ such
that $r$ induces a bijection of $\wtI_{c,1}$ with $I_c$. 
Putting $\wtI_1=\wtI_{nc,1}\cup\wtI_{c,1}$, we
get an identification
\[
T_{L/\CC}=\bigoplus_{\tilde\iota\in\wtI_1}\Gm{/\CC}.
\]
Under this identification the factor $\Gm{/\CC}$ corresponding to
$\tilde\iota$ acts on the factor of $V_{L/\CC}$ corresponding to
$\tilde\iota$ by multiplication in case $\tilde\iota\in\wtI_1$ and by
multiplication with the inverse if $\tilde\iota\not\in\wtI_1$.

The quasi-cocharacter 
$\nu_L\colon\Gm{/\CC}\rightarrow T_{L/\CC}$ is defined by 
\[
\nu_L(z)_{\tilde\iota}=
\left\{\begin{array}{ll} 
1&\textrm{if }\tilde\iota\in\wtI_{c,1}\\
z^{-1/4}&\textrm{if }\tilde\iota\in\wtI_{nc,1}.
\end{array}\right.
\]
We put $\wtV=W\otimes_L V_L$ as representation of $\wtG^\der\times T$
and again let $\wtG$ be the image of $\wtG^\der\times T$ in $\GL(\wtV)$. 
As before, the derived group of $\wtG$ is $\wtG^\der$ and
the projection of $(\nu,\nu_L)$ from $(\wtG^\der\times T)_{/\CC}$ to
$\wtG_{/\CC}$ is a true cocharacter $\tmu$ of $\wtG_{/\CC}$. 
We obviously have equalities
$\tilde\pi^\ad\circ\tmu=\tilde\pi^\ad\circ\nu=\pi^\ad\circ\mu$, so  
this settles the construction of the couple $(\wtG,\tmu)$ in this
case. 

The preceding discussion establishes the following theorem. 
\begin{thm}\label{cstr_wtG-tmu}
Let $A/\CC$ be a simple abelian variety, $(G_A,h_A)$ its
Mumford--Tate datum and $\mu_A$ the Hodge cocharacter
associated to $A$. 
Assume that $A$ is essentially Mumford--Tate unliftable
and that $G_A$ is of type $D_k^\HH$ with $k\geq4$.
Then there exist an algebraic group $\wtG$ with $\wtG^\der$ simply
connected, an identification $\wtG^\ad=G_A^\ad$ and a cocharacter
$\tmu\colon\Gm{/\CC}\rightarrow\wtG_{/\CC}$ such that
$\tilde\pi^\ad\circ\tmu=\pi^\ad\circ\mu_A$. 
\end{thm}
\subsection{Remark.}\label{torus_compact}
In all cases above, the group $T_{L/\RR}$ occurring in the construction of
$\wtG$ is compact. 
It thus follows that $\wtG^\ab_{/\RR}$ is compact. 
It is left to the reader to construct an isomorphism
$\wtG^\ab\cong T_L$ and to compute the composite
$T_L\subset\wtG\rightarrow\wtG^\ab\cong T_L$.
\subsection{}\label{LADH_intro}
Suppose that $F\subset\CC$ is a number field and that $A$ is an abelian
variety over $F$ with connected Mumford--Tate group such that 
the conditions of the theorem are verified for $A_{/\CC}$ and
its Mumford--Tate datum $(G_A,h_A)$. 
Let $\wtG$ and $\tmu$ be as in the conclusion of the theorem. 

Let $\wtG^\ab=\wtG/\wtG^\der$ and 
define $\mu^\ab$ to be the composite of the cocharacter 
$\tmu\colon\Gm{/\CC}\rightarrow\wtG_{/\CC}$ with the natural projection
$\wtG_{/\CC}\rightarrow\wtG^\ab_{/\CC}$. 
Define
\begin{align*}
h^\ab\colon S&\rightarrow\wtG^\ab_{/\RR}\\
z&\mapsto\mu^\ab(z)\overline{\mu^\ab(\bar z)}.
\end{align*}
The remark~\ref{torus_compact} implies that $\wtG^\ab_{/\RR}$ is
compact.
As the weight $w_h\colon\Gm{/\RR}\rightarrow\wtG^\ab_{/\RR}$ is
trivial, \cite[1.1.18(b)]{deligne:shimura2} implies
that $h^\ab$ defines a polarizable Hodge structure on any $\QQ$-linear
representation $V^\ab$ of $\wtG^\ab$. 

Fix a faithful rational representation $V^\ab$ of $\wtG^\ab$.  
It follows from proposition~A.1 of \cite{deligne:motifs_taniyama} and
the remark preceding it that the Hodge structure on $V^\ab$ is the
Betti realization of an absolute Hodge motive $M^\ab_{/\bQ}$ belonging
to $\CM_{\bQ}$. 
Note that using the isomorphism $\wtG^\ab\cong T_L$ from
remark~\ref{torus_compact}, this motive can also be constructed
explicitly. 
Replacing $F$ by a finite extension and fixing an embedding 
$F\subset\bQ$ we can assume that $M^\ab_{/\bQ}$ descends
to a motive $M^\ab$ over $F$. 
At the cost of a further finite extension of $F$, 
we can also assume that the Mumford--Tate group of
$M^\ab$ is connected and hence contained in $\wtG^\ab$. 
This implies that for every $p$, the $p$-adic realization of $M^\ab$
factors through a map $\rho_p^\ab\colon\cG_F\rightarrow\wtG^\ab(\Qp)$. 
Note that the $\rho_p^\ab$ form a compatible system
of representations (see~\ref{frobenii}) 
and that the $\rho_p^\ab$ do not depend on the
choice of the representation $V^\ab$ of $\wtG^\ab$. 

For each prime number $p$, let $\Sigma_p$ be the set of $p$-adic
places of $F$. 
There is a finite set $\Sigma$ of places $v$ of $F$ such that $A$ and
$M^\ab$ have good reduction outside $\Sigma$.
For every prime number $p$, the representations $\rho_{A,p}$ 
and $\trho_p^\ab$ are unramified at all $v\not\in\Sigma\cup\Sigma_p$. 
For any finite extension $F'$ of $F$, let $\Sigma'_p$ be the set of
$p$-adic places of $F'$ and let $\Sigma'$ be the set of places lying
over the $v\in\Sigma$. 

Let $G'=G_A^\ad\times\wtG^\ab$ and let $\pi'\colon\wtG\rightarrow G'$ be
defined by the projections of $\wtG$ onto $\wtG^\ab$ and
$G_A^\ad=\wtG^\ad$. 
For each $p$, define
$\rho_p^\ad=\pi^\ad\circ\rho_{A,p}\colon\cG_F\rightarrow G_A^\ad(\Qp)$
and put 
\[
\rho_p'=(\rho_p^\ad,\rho_p^\ab)\colon\cG_F\rightarrow G'(\Qp). 
\]
\begin{cor}\label{constr_LADH}
There exist a finite extension $F'$ of $F$ and a system of geometric
Galois representations $\trho_p\colon\cG_{F'}\rightarrow\wtG(\Qp)$
lifting the restrictions $\rho_p'\colon\cG_{F'}\rightarrow G'(\Qp)$.
The system $(\trho_p)$ can be chosen in such a way that each
$\trho_p$ is unramified at $v$ for each $v\not\in\Sigma'\cup\Sigma_p'$.

For each prime number $p$, the representation $\trho_p$ is a weak
geometric lift of $\rho_{A,p}$ and it is essentially geometrically
unliftable. 
\end{cor}
\begin{proof}
Since $G_A$ is the Mumford--Tate group of $A$, the adjoint
representation of $G_A$ gives rise to a motive 
$M^\ad$ belonging to the subcategory of $\Mot_\AH(F)$ which is
$\otimes$-gene\-ra\-ted by $h^1(A)$ and the Tate motive. 
The Mumford--Tate group of $M^\ad$ is the adjoint group
$G_A^\ad=\wtG^\ad$. 

Consider the object $M'=M^\ad\times M^\ab$ of $\Mot_\AH(F)$. 
Its Mumford--Tate group is contained in $G'$ and the
Hodge cocharacter $\mu'\colon\Gm{/\CC}\rightarrow G'_{/\CC}$
associated to its Betti realization is the product
$(\mu^\ad,\mu^\ab)$. 
It follows that $\mu'=\pi'_{/\CC}\circ\tmu$. 
On the other hand, for each $p$, the Galois representation
on the $p$-adic realization of the motive $M'$ is the above map 
$\rho_p'$. 
The motive $M'$ is an abelian motive, so it follows from
\cite[Theorem~0.3]{blasius:p-adic_hodge} that the $p$-adic comparison
maps linking the $p$-adic and the DeRham realizations of $M'$ are
compatible with absolute Hodge classes. 
The main theorem~2.1.7 of \cite{wintenberger:relevement} implies that
there is a finite extension $F'\supset F$ such that every 
$(\rho_p')_{|\cG_{F'}}$ lifts to a representation 
$\trho_p\colon\cG_{F'}\rightarrow\wtG(\Qp)$ satisfying the conditions
of the corollary.

It is clear from the above that, for each $p$, the map $\trho_p$
provides a weak geometric lift of $\rho_{A,p}$. 
As $\wtG^\der$ is simply connected, it is also clear that the
$\trho_p$ are essentially geometrically unliftable. 
\end{proof}
\begin{dfn}\label{def_lifted_abelian_DkH}
The representations $\trho_p$ constructed as in the corollary from
Mum\-ford--\-Tate decomposed abelian varieties will be
called \emph{representations of lifted abelian $D_k^\HH$-type}. 
\end{dfn}
\begin{cor}\label{simply_connected_wgls_exist}
Let $F\subset\CC$ be a number field and $A/F$ an abelian variety with
connected Mumford--Tate group $G_A$ and associated system of Galois
representations $(\rho_{A,p})$.
Then there exist an linear algebraic group $\wtG$ over $\QQ$ such that
$\wtG^\der$ is the universal cover of $G_A^\der$, 
a finite extension $F'$ of $F$
and a system of weak geometric lifts 
\[\trho_p\colon\cG_{F'}\rightarrow\wtG(\Qp)\] of the restrictions to
$\cG_{F'}$ of the $\rho_{A,p}$. 
\end{cor}
\begin{proof}
By theorem~\ref{unliftable_lift_exists} and
proposition~\ref{descent_MT-adjoint}, there exist a number field
$F'$ and an essentially M-T unliftable and M-T decomposed weak
Mumford--Tate lift $B/F'$ of $A$ with connected Mumford--Tate group $G_B$.
It follows from proposition~\ref{weak_MT_lift-galois_lift} that the
system of Galois representations 
$(\rho_{B,p})$ associated to $B$ is a system of weak geometric lifts
of the system $(\rho_{A,p})$. 

Let $B_{/\CC}\sim\prod B_{i/\CC}$ be the isogeny decomposition as in
definition~\ref{def_MT-reduced}. 
After enlarging $F'$, we can assume that this decomposition exists over
$F'$. 
In that case, each $B_i$ is Mumford--Tate decomposed with
$G_{B_i}^\ad$ simple and $G_B^\ad$ is the product of the
$G_{B_i}^\ad$. 
For each $i$, let 
$\rho_{B_i,p}\colon\cG_{F'}\rightarrow G_{B_i}(\Qp)$ be
the $p$-adic Galois representation associated to $B_i$. 
For every $i$ such that $G_{B_i}$ is not of type $D_k^\HH$, put
$\wtG_i=G_{B_i}$ and $\trho_{i,p}=\rho_{B_i,p}$. 
For each $i$ such that $G_{B_i}$ is of type $D_k^\HH$ let
$\trho_{i,p}\colon\cG_{F'}\rightarrow\wtG_i(\Qp)$ be the
representation of lifted abelian $D_k^\HH$-type resulting from 
corollary~\ref{constr_LADH}, replacing $F'$ by a finite extension
again if necessary. 
Let $\wtG=\prod\wtG_i$ and let
$\trho_p\colon\cG_{F'}\rightarrow\wtG(\Qp)$ be the product of the maps
$\trho_{i,p}$.

For each prime number $p$, the representation $\trho_p$ is a weak
geometric lift of $\rho_{B,p}$ and therefore also of $\rho_{A,p}$. 
As $\wtG^\der$ is simply connected, the corollary follows. 
\end{proof}
\begin{cor}\label{LADHs_dominate_all}
Let $F$ be a number field, $M$ an abelian motive and 
$\rho'\colon\cG_F\rightarrow G'(\Qp)$ a weak geometric lift of the
associated $p$-adic Galois representation. 
Then, after replacing $F$ by a finite extension and for any linear
representation $V'_p$ of $G'_{/\Qp}$, the 
representation of $\cG_F$ on $V'_p$ deduced from $\rho'_p$ lies in 
the tannakian subcategory of 
$\Rep_{\Qp}(\cG_F)$ generated by the $p$-adic representations
associated to abelian motives and the representations of lifted
abelian $D_k^\HH$-type.
\end{cor}
\begin{proof}
It is 
sufficient to prove the corollary for one fixed faithful self-dual
representation $V'_p$ of $G'$. 

Let $A$ be the essentially Mumford--Tate unliftable and
Mumford--Tate decomposed abelian variety supplied by
applying corollary~\ref{AM_to_MTud-AV} to $M$, so $M$ is an
object of $\langle h^1(A),\QQ(1)\rangle$.
This inclusion corresponds to a map $G_A\rightarrow G_M$.
Suppressing isogeny factors of $A$, we can assume that the induced map
$G_A^\der\rightarrow G_M^\der$ is an isogeny. 
After replacing $F$ by a finite extension, define the group $\wtG$ and
the representation $\trho_p\colon\cG_F\rightarrow\wtG(\Qp)$ as in the
proof of~\ref{simply_connected_wgls_exist}. 

As $\wtG^\der$ is the simply connected cover of $G_A^\der$, the
isogeny 
\[
\wtG^\der\rightarrow G_A^\der\rightarrow G_M^\der
\] 
lifts to an isogeny $\wtG^\der\rightarrow G^{\prime\der}$.
This results in an identification of the adjoint groups
$G^{\prime\ad}=G_M^\ad=G_A^\ad=\wtG^\ad$ and, by construction, the
projections to $G_M^\ad(\Qp)$ of the representations $\rho'_p$,
$\rho_{M,p}$, $\rho_{A,p}$ and $\trho_p$ all coincide.
The corollary now follows from
lemma~\ref{if_wgl_then_tannakian} 
with $G_1=\wtG$, $V_1$ any faithful self-dual representation of $G_1$,
$G_2=G'$ and $V_2=V'_p$.
\end{proof}
\section{Galois representations of lifted abelian $D_k^\HH$-type}
\label{sect_LADH}
\subsection{}\label{intro_LADH}
We revert to the notations of~\ref{LADH_intro}, in particular $A/F$ is
essentially 
Mumford--Tate unliftable, $G_A$ is of type $D_k^\HH$ and $\wtG$ is the
group constructed in theorem~\ref{cstr_wtG-tmu}.
We will also assume $A$ to be Mumford--Tate decomposed. 
For the first part of this section, we fix a prime number $p$ and
assume that there exists a weak geometric lift
$\trho_p\colon\cG_F\rightarrow\wtG(\Qp)$ of $\rho_{A,p}$ with $\wtG$
as in theorem~\ref{cstr_wtG-tmu}. 
It follows from corollary~\ref{constr_LADH} that such a lifting exists
after replacing $F$ by a finite extension and that, for $F$ big
enough, $\trho_p$ is part of a system of $p$-adic representations for
variable $p$, but we do not assume for the moment that such a system
exists. 
\begin{thm}\label{LADH_non-abelian}
Let $\trho_p\colon\cG_F\rightarrow\wtG(\Qp)$ be a representation of
lifted abelian $D_k^\HH$-type, constructed as weak geometric lift of
the representation $\rho_{A,p}$ as described in~\ref{intro_LADH}, with
$k\geq5$. 
Assume moreover that $\rho_{A,p}(\cG_F)\subset G_A(\Qp)$ is Zariski
dense. 

Let $\wtV_p$ be a faithful $\Qp$-linear representation of
$\wtG_{/\Qp}$. 
Then there is no finite extension $F'$ of $F$ such that the
representation of $\cG_{F'}$ on $\wtV_p$ induced by $\trho_p$ belongs to the
category $\RepAV{\Qp}{(\cG_{F'})}$. 
\end{thm}
\begin{proof}
Replace $F$ by a finite extension and assume that the representation
of $\cG_F$ on $\wtV_p$ belongs to $\RepAV{\Qp}{(\cG_{F})}$. 
After further enlarging $F$ and applying
proposition~\ref{weak_MT_lift-galois_lift} and 
remark~\ref{B'_in_thm_exists}, it follows
that there is an essentially Mumford--Tate unliftable and
Mumford--Tate decomposed abelian variety $B/F$ such that $\wtV_p$
belongs to $\langle\rH_\et^1(B_{\bar F},\Qp),\RepCM{}{}\rangle$. 

Let $B\sim\prod_{j=1}^{m-1} B_j$ be the decomposition from the
definition~\ref{def_MT-reduced}, so that
$G_B^\der=\prod G_{B_j}^\der$, each $G_{B_j}^\ad$ is simple and
$G_B^\ad=\prod G_{B_j}^\ad$.
For $j=1,\cdots,m-1$, let $\rho_j=\rho_{B_j,p}$, put
$V_j=\rH_\et^1(B_{j/\bar F},\Qp)$ and let $H_j$ be the
Zariski closure of $\rho_j(\cG_F)$ in $G_{B_j}$. 
Let $V_m$ be an object of $\RepCM{}{}$ such that $\wtV_p$
belongs to $\langle V_1,\ldots,V_m\rangle$, let $H_m$ be the Zariski
closure of the image of the Galois representation on $V_m$ and let
$\rho_m\colon\cG_F\rightarrow H_m(\Qp)$ be the corresponding morphism.
Finally, let $H\subset\prod_{j=1}^m H_j$ be the Zariski closure of
the image of 
\[
(\rho_1,\ldots,\rho_m)\colon\cG_F\rightarrow \prod_{j=1}^m H_j(\Qp)
\] 
and $\sigma\colon\cG_F\rightarrow H(\Qp)$ the induced representation. 
At the cost of a further finite extension of $F$, we can assume $H$ to
be connected. 
The fact that $\wtV_p$ belongs to the tannakian category 
generated by the $V_j$ implies that there is a surjection
$\pi_p\colon H\rightarrow\wtG_{/\Qp}$ such that
$\trho_p=\pi_p\circ\sigma$. 

For any $p$-adic place $\bar v$ of $\bar F$, the cocharacter
$\mu_{\trho_p,\bar v}\colon\Gm{/\Cp}\rightarrow\wtG_{/\Cp}$
associated to the Hodge--Tate decomposition corresponding to
$\left(\trho_p\right)_{|\cI_{\bar v}}$ lifts to the cocharacter
$\mu_{\sigma,\bar v}\colon\Gm{/\Cp}\rightarrow H_{/\Cp}$ associated to
the Hodge--Tate decomposition corresponding to $\sigma_{|\cI_{\bar v}}$. 
There is at least one simple factor $G'_{/\Cp}$ of $\wtG^\ad_{/\Cp}$
to which $\mu_{\trho_p,\bar v}$ projects non-trivially.
It follows from  \cite[Proposition~7]{wintenberger:motifs} and the
fact that \cite[Theorem~0.3]{blasius:p-adic_hodge}
implies \cite[Conjecture~1]{wintenberger:motifs} that this projection is
dual to the root $\alpha_k$ of this factor. 

Let $H'_{/\Cp}$ be the simple isogeny factor of $H^\der_{/\Cp}$ which
surjects onto $G'_{/\Cp}$, such a factor exists by hypothesis. 
As $\wtG^\der$ is simply connected, $H'_{/\Cp}$ is simply connected as
well. 
The Hodge--Tate cocharacter $\mu_{\sigma,\bar v}$ lifts to a
quasi-cocharacter  $\mu'_{\bar v}$ of $H'_{/\Cp}$.
This quasi-cocharacter is still dual to the vertex $\alpha_k$ of the
Dynkin diagram. 
Over $\Cp$, any faithful representation $W'$ of $H'_{/\Cp}$ whose
highest weights 
are fundamental weights contains a
direct factor with highest weight $\varpi_{k-1}$ or $\varpi_k$ and it
follows from lemma~\ref{weights_nu} applied to $H'_{/\Cp}$, with
$\mu'_{\bar v}$ playing the role of $\nu$, that $\mu'_{\bar v}$ has at
least three weights on this factor. 
We will show that this leads to a contradiction. 

As $W=\bigoplus_{j=1}^{m-1}V_j$ is the Galois representation on
$\rH_\et^1(B_{\bar F},\Qp)$, the Hodge--Tate cocharacter
$\mu_{\sigma,\bar v}$ acts on $W\otimes\Cp$ with two weights. 
By construction, $W$ is a faithful representation of $H^\der$ so
$W\otimes_{\Qp}\Cp$ is a faithful representation of $H'_{/\Cp}$.
Let $W'$ be the direct sum of the direct factors of the representation
of $H'_{/\Cp}$ on $W\otimes\Cp$ on which $H'_{/\Cp}$ acts non-trivially. 
The quasi-cocharacter $\mu'_{\bar v}$ then acts with exactly two weights on
$W'$. 
It follows from \cite[1.3.7]{deligne:shimura2} that for every
irreducible direct factor $W''$ of $W'$, the highest
weight is a fundamental weight of $H'_{/\Cp}$. 
We have shown above that this implies that $\mu'_{\bar v}$ has at least three
weights on $W''$, so we arrive at the contradiction we were looking
for. 
\end{proof}
\subsection{Remark.}
It follows from theorems~\ref{MT-group_unliftable-AV}
and~\ref{unliftable_lift_exists} and their proof, that for every
$k\geq4$ there exists an essentially M-T unliftable and M-T decomposed
abelian variety $A/\CC$ such that $G_A$ is of type $D_k^\HH$. 
It follows from \cite[Theorem~1.7]{noot:av_etc} that there also exists such
an abelian variety which can be defined over a number field $F$ and for
which $\rho_{A,p}(\cG_F)\subset G_A(\Qp)$ is Zariski dense. 
Thus, for every $k\geq5$, there is an abelian variety for which the
hypothesis of the theorem is fulfilled. 
\subsection{The case $k=4$.}
We consider the case of Galois representations of lifted abelian
$D_4^\HH$-type, so we keep all
the notations of~\ref{intro_LADH} but fix $k=4$. 
In this case, the method of
the proof of the above theorem gives a more limited statement. 
The reason for this lies in the facts that the representation of
$\cG_F$ on $\wtV_p$ may be reducible and that the group
$\wtG_{/\CC}$ does have faithful representations in which the Hodge
cocharacter acts with only two weights. 

It is left to the reader to verify the following statement, whose
proof is completely analogous to the proof of
theorem~\ref{LADH_non-abelian}. 
\begin{prop}
Let notations and hypotheses be as in theorem~\ref{LADH_non-abelian},
but with $k=4$.  

Let $\wtV$ be a faithful $\QQ$-linear representation of $\wtG$.
Then there is no finite extension $F'$ of $F$ such that the
representation of $\cG_{F'}$ on $\wtV\otimes\Qp$ induced by $\trho_p$
belongs to the category $\RepAV{\Qp}{(\cG_{F'})}$. 
\end{prop}
\subsection{Remark.}
In the situation of the proposition, assume that $p$ is a prime number
such that the image of $\cG_{\Qp}$ in the automorphism group of the
Dynkin diagram is equal to the image of $\cG_\QQ$ in this automorphism
group. 
Then the statement of the proposition (and of
theorem~\ref{LADH_non-abelian}) holds for 
every $\Qp$-linear representation $\wtV_p$ of $\wtG_{/\Qp}$. 
If $\cG_\QQ$ acts on the Dynkin diagram of $G_{\bQ}$ through a cyclic
group, then the images of $\cG_{\QQ}$ and $\cG_{\Qp}$ in the
automorphism group of the Dynkin diagram coincide for infinitely many
$p$.

On the other hand, there also exist examples of Galois representations
of lifted abelian $D_4^\HH$-type where the proof of the theorem fails
for every $p$. 
One can find an example of this situation by taking an abelian variety
for which $\cG_\QQ$ acts on the Dynkin diagram of the Mumford--Tate
group through $(\ZZ/2\ZZ)^2$. 
For every prime number $p$, the image of $\cG_{\Qp}$ in the
automorphism group of the Dynkin diagram is then of order at most
$2$. 
\subsection{Frobenius elements.}\label{frobenii}
Let $v$ be a valuation of $F$. 
By $\Fr_v\in\cG_F$, we will denote a geometric Frobenius element at $v$, 
i.\ e.\ an element of the decomposition group of a place $\bar v$ of
$\bar F$ lying over $v$ such that, on the residue field $\bar k_v$ of
$\bar F$ at $\bar v$,
$\Fr_v$ induces the \textit{inverse} of the map $x\mapsto x^{q_v}$, 
where $q_v$ is the order of the residue field $k_v$ of $F$ a $v$.
Note that $\Fr_v$ is defined only up to conjugation and up to
multiplication by an element of $\cI_{F,\bar v}$.
This implies that the image of $\Fr_v$ in a representation which is
unramified at $v$ is defined up to conjugation and that 
the eigenvalues and the characteristic polynomial of the image of
$\Fr_v$ in such a representation are well defined. 

Let $V_p$ be a $\Qp$-linear representation of $\cG_F$ on an \'etale
cohomology group of a proper and smooth $F$-variety and let $v$ be a
place where this variety has good reduction. 
It follows from the Weil conjectures, proved by Deligne, that the
eigenvalues of a Frobenius element $\Fr_v$ at $v$ are algebraic
integers and that all complex absolute values of all these eigenvalues
coincide. 
If the Fontaine--Mazur conjecture is true, then every geometric
representation of $\cG_F$ should have this property. 
We will verify this for the representations of lifted abelian
$D_k^\HH$-type. 
It is also shown that $\Fr_v$ acts semi-simply in any representation
of lifted abelian $D_k^\HH$-type, a property which is conjectured for
all representations coming from the cohomology of algebraic
varieties. 

The Weil conjectures also imply that, for varying $p$, the $p$-adic
\'etale cohomology groups of a proper and smooth variety form a
\emph{compatible system} $(V_p)$ of Galois representations.  
This means that there is a finite set $\Sigma$ of valuations of $F$
such that $V_p$ is unramified at $v$ for all $v\not\in\Sigma$ with
$v(p)=0$ and that for all $v\not\in\Sigma$ there is
a polynomial $P^v\in\QQ[X]$ which is equal to the characteristic
polynomial of $\Fr_v$ acting on $V_p$ for all $p$ with $v(p)=0$. 
We will prove a result in this direction for the representations of
lifted abelian $D_k^\HH$-type. 

It should be pointed out that this property does not follow from the
conjecture of Fontaine and Mazur. 
However, combined with the Mumford--Tate conjecture, the Fontaine--Mazur
conjecture implies that
for any geometric $p$-adic representation $W_p$ of $\cG_F$, there should exist
a number field $E$ and a compatible system (indexed by the primes
$\mathfrak{p}$ of $E$) of $E_{\mathfrak{p}}$-linear
representations of $\cG_F$ such that $W_p$ occurs in this system. 
The characteristic polynomials of $\Fr_v$ acting on the $V_p$ would
then lie in $E[X]$ and be independent of $p$, 
for all $p$ with $v(p)=0$.
\begin{prop}\label{weight_LADH}
Let $F\subset\CC$ be a number field and let 
$\trho_p\colon\cG_F\rightarrow\wtG(\Qp)$ be a
representation of lifted abelian $D_k^\HH$-type, with $k\geq4$.
Let $\Sigma$ be the set of places of bad reduction defined
in~\ref{LADH_intro} and 
let $\wtV$ be an irreducible $\QQ$-linear representation of $\wtG$. 

Then, for each $v\not\in\Sigma\cup\Sigma_p$, the Frobenius element
$\trho_p(\Fr_v)$ acts semi-simply on $\wtV_p=\wtV\otimes\Qp$ and 
its eigenvalues are algebraic integers with all complex absolute
values equal to $1$. 
\end{prop}
\begin{proof}
It suffices to prove the proposition for the faithful representation $\wtV$
constructed in~\ref{case_k_even} and~\ref{case_k_odd} (according to
the parity of $k$) and after replacing $F$ by a finite extension. 
Let $\trho_p$ be lifted from the representation associated to the
M-T decomposed and M-T unliftable abelian variety
$A/F$, with Mumford--Tate group 
$G_A$. 
The derived group $G_A^\der$ is of the form $\Res_{K_0/\QQ}G^{s,\der}$
and we have 
$\wtG^\der=\Res_{K_0/\QQ}\wtG^{s,\der}$ where $\wtG^{s,\der}$ is the
universal cover of $G^{s,\der}$.
By construction, the representation $\wtV$ carries a structure of
$K_0$-vector space compatible with the structure of a Weil restriction
on $\wtG^\der$, so the action of $\wtG$ on $\wtV$ commutes with the
action of $K_0$. 

It is enough to prove the proposition for the representation
of $\cG_F$ on the tensor product $\wtW_p=\wtV_p\otimes_{K_{0,p}}\wtV_p$, 
where $K_{0,p}=K_0\otimes_\QQ\Qp$. 
The representation of $\wtG^\der_{/\Qp}$ on $\wtW_p$ factors
through a representation of $G_{A/\Qp}^\der$.
In fact, $\wtW_p$ is
a faithful representation of $G_{A/\Qp}^\der$, but this plays no role
in what follows.  
By 
theorem~\ref{thm_non_spinorial_lift_from_AM}, $\wtW_p$ occurs in
the $p$-adic realization of an abelian motive $M$. 
The proof of~\ref{thm_non_spinorial_lift_from_AM} even shows that we
can find such a motive $M$ in the subcategory of $\AV_F$ generated by
$h^1(A)$ and $\CM_F$ so it follows that $M$ has potentially good
reduction at $v$ for each $v\not\in\Sigma$. 
After replacing $F$ by a finite extension, we can assume that $M$ has
good reduction at all $v\not\in\Sigma$. 
This implies that for any valuation 
$v\not\in\Sigma\cup\Sigma_p$, the Frobenius element
$\trho_p(\Fr_v)$ acts semi-simply on $\wtW_p$ and that each eigenvalue
$\lambda$ is an algebraic integer with all complex absolute
values of the form $\NN v^{w(\lambda)/2}$, where $\NN v$ is the
cardinal of the residue field of $F$ at $v$ and $w(\lambda)$ is an  
integer. 

To show that all $w(\lambda)$ are equal to $0$, it suffices to show that
the Betti realization of $M$ is pure of weight $0$. 
This realization can be determined as follows. 
The cocharacter $\tmu$ constructed in~\ref{cstr_wtG-tmu} defines a map
$\tih_\CC\colon S_{/\CC}\rightarrow\wtG_{/\CC}$ given by 
$\tih_\CC(z,\bar z)=\tmu(z)\overline{\tmu(\bar z)}$. 
This map descends to a map $\tih\colon S\rightarrow\wtG_{/\RR}$ lifting
$h_A^\ad\colon S\rightarrow G_A^\ad$ and 
these data define a Hodge structure on $\wtV$ on which $K_0$ acts by
endomorphisms. 
The tensor product $\wtW=\wtV\otimes_{K_0}\wtV$ is a 
representation of $G_A^\der$, and the Hodge structure on $\wtW$ derived
from the
Hodge structure on $\wtV$ is the Hodge structure on the Betti
realization of the abelian motive $M$. 
Since $\tmu$ and its complex conjugate are inverse to each other, it
follows that the Hodge structure on $\wtV$ is of weight $0$ and thus
the same thing holds for $\wtW$.
\end{proof}
\subsection{}\label{intro_char_pol}
We keep the notations of proposition~\ref{weight_LADH}.
In particular,
$\Sigma$ is the set of places of bad reduction defined
in~\ref{LADH_intro}. 
For the representation $\wtV$ however, we fix the $\QQ$-linear
representation of $\wtG$ 
constructed in~\ref{case_k_even} or~\ref{case_k_odd} as in the
proof of proposition~\ref{weight_LADH}.
It follows from corollary~\ref{constr_LADH} that, after replacing $F$
by a finite extension, there is a system of weak geometric lifts
$\trho_p\colon\cG_F\rightarrow\wtG(\Qp)$ of the $\rho_{A,p}$, for
varying $p$. 
From now on we will assume that we dispose of such a system. 

The rest of the paper concerns the variation of the characteristic
polynomials of the $\trho_p(\Fr_v)$ on the $\wtV_p=\wtV\otimes\Qp$, for
a fixed valuation $v\not\in\Sigma$ of $F$ and
varying $p$ with $v(p)=0$.  
We will deduce our result on the characteristic polynomials
from the statement~\ref{classes_Fr} which is independent of the choice
of a representation of $\wtG$. 

As in the proof of proposition~\ref{weight_LADH}, assume that
$\trho_p$ is lifted 
from the representation associated to the Mumford--Tate decomposed
and essentially M-T unliftable abelian variety $A/F$, with Mumford--Tate
group $G_A$. 
Without loss of generality, we can assume that $G_A$ is of the type
constructed in remark~\ref{better_unliftable_lift}. 
We also keep the notations $G_A^\der=\Res_{K_0/\QQ}G^{s,\der}$ and
$\wtG^\der=\Res_{K_0/\QQ}\wtG^{s,\der}$, for a totally real number
field $K_0$ and semi-simple groups $G^{s,\der}$ and $\wtG^{s,\der}$
over $K_0$. 
As explained in~\ref{better_unliftable_lift}, there is a totally
imaginary quadratic extension $L$ of $K_0$ such that $G_A$ is
isogenous to a subgroup of $G_A^\der\times L^\times$. 
It follows from loc.\ cit.\ that $L$ is contained in the centre of
$\End^0(A_{/\CC})=\End(A_{/\CC})\otimes_\ZZ\QQ$. 
In this case, there exists a linear algebraic group $\wtG^s$ over
$K_0$ such that $\wtG=\Res_{K_0/\QQ}\wtG^s$

Let $\wtN\subset\wtG^\der\subset\wtG$ be the centre of $\wtG^\der$. 
It is the kernel of the natural map
$\wtG\rightarrow G'=\wtG^\ad\times\wtG^\ab$.
There is an isomorphism $\wtN\cong\Res_{K_0/\QQ}\wtN^s$, where
$\wtN^s$ is the centre of $\wtG^{s,\der}$, a finite group scheme over
$K_0$ of order $4$, geometrically isomorphic to $(\ZZ/2\ZZ)^2$ if $k$
is even and to $\ZZ/4\ZZ$ if $k$ is odd. 
Note that, by construction, each $\trho_p$ is determined up to a
finite character $\cG_F\rightarrow\wtN(\Qp)$ which is unramified outside
$\Sigma\cup\Sigma_p$. 
This ambiguity explains the fact that in the next proposition the
characteristic polynomials may vary with $p$. 

For each valuation $v$ of $F$, each prime number $p$ and each
$\epsilon\in\wtN(\Qp)$, let
$\wtP^v_{\epsilon,p}(X)\in\Qp[X]$ be the characteristic polynomial 
\[
\wtP^v_{\epsilon,p}(X)=
\det\nolimits_{\wtV_p}\left(\epsilon\trho_p(\Fr_v)-X\cdot\id\right).
\]
We write $\wtP^v_p(X)=\wtP^v_{1,p}(X)$ for the characteristic polynomial of
$\trho_p(\Fr_v)$. 
\begin{prop}\label{char_pols_Fr}
Assume that for some, hence any, prime number $\ell$, the rank of the
Zariski closure of the image of $\rho_{A,\ell}$ is equal to the rank
of $G_A$. 
Then there exist a set $\Sigma_{\Fr}\supset\Sigma$ of valuations of $F$ of
Dirichlet density $0$ and, 
for each $v\not\in\Sigma_{\Fr}$, a polynomial $\wtP^v\in\QQ[X]$ such
that for every prime number $p$ with $v(p)=0$, one has
$\wtP^v_{\epsilon,p}=\wtP^v$ 
for some $\epsilon=\epsilon_p\in\wtN(\Qp)$.  
\end{prop}
\begin{cor}
Let notations and hypotheses be as in the proposition. 
Let $M^\mathrm{nor}$ be a normal closure of $K_0/\QQ$.
For each prime number $p$, let $M^\mathrm{nor}_p$ be the image of
$M^\mathrm{nor}$ in $\bQp$ and put $M_p=M^\mathrm{nor}_p\cap\Qp$. 
Then, for every $v\not\in\Sigma_{\Fr}$, the
characteristic polynomial $\wtP^v_p$ of
$\trho_p(\Fr_v)$ lies in $M_p[X]$. 
\end{cor}
\begin{proof}[Proofs]
The $\ell$-independence of the rank of the Zariski closure
of the image of $\rho_{A,\ell}$ follows from
\cite[2.2.4]{serre:resume84-85}. 
The corollary easily follows from the proposition. 
To prove the proposition we will reformulate the result in terms of
the the quotient variety of $\wtG$ by a subgroup
$\Aut'\subset\Aut(\wtG^\der)$.

If $k\geq5$, we put $\Aut'(\wtG^{s,\der})=\Aut(\wtG^{s,\der})$ and
$\Out'(\wtG^{s,\der})=\Out(\wtG^{s,\der})$. 
If $k=4$ then $\Out(\wtG^{s,\der})\cong\mathrm{S}_3$, acting naturally
on the
vertices $\alpha_1,\alpha_3,\alpha_4$ of the Dynkin diagram.
Let $\Out'\subset\Out(\wtG^{s,\der})$ be the stabiliser of $\alpha_1$
and $\Aut'(\wtG^{s,\der})$ the inverse image of $\Out'$ in
$\Aut(\wtG^{s,\der})$. 

In either case, $\Aut'(\wtG^{s,\der})$ is an extension of
$\Out'(\wtG^{s,\der})\cong\ZZ/2\ZZ$ by $\wtG^{s,\ad}$. 
It acts on the centre of $\wtG^{s,\der}$ through its quotient
$\Out'(\wtG^{s,\der})$ and by 
going through the constructions in~\ref{constr_tildes}
to~\ref{case_k_odd} (sub)case by (sub)case, it is not very difficult
to check that this action extends to an action of $\Out'(\wtG^{s,\der})$
on the centre of $\wtG^s$. 
We thus obtain an action of $\Aut'(\wtG^{s,\der})$ on $\wtG^s$ and,
taking Weil restrictions, an action of
$\Aut'=\Res_{K_0/\QQ}\Aut'(\wtG^{s,\der})$ on $\wtG$. 

Let $\Cl(\wtG)$ be the categorical quotient of $\wtG$ by this
action, see \cite[Chap\-ter~1]{mumfog:GIT}.
This means that if $R=\Gamma(\wtG,\mathcal{O}_{\wtG})$ is the
affine coordinate ring of $\wtG$ then the quotient is given by
$\Cl(\wtG)=\Spec(R^{\Aut'})$.
For each prime number $p$, we denote by
$\Cl(\trho_p)\colon\cG_F\rightarrow\Cl(\wtG)(\Qp)$ the map deduced
from $\trho_p$. 
The proposition~\ref{char_pols_Fr} then follows from the following,
more precise, statement. 
\end{proof}
\begin{prop}\label{classes_Fr}
Assume that for some, hence any, prime number $\ell$, the rank of the
Zariski closure of the image of $\rho_{A,\ell}$ is equal to the rank
of $G_A$. 
Then there are a set $\Sigma_{\Fr}\supset\Sigma$ of valuations of $F$
of Dirichlet density $0$ and elements $\Cl(\wtFr_v)\in\Cl(\wtG)(\QQ)$,
for $v\not\in\Sigma_{\Fr}$, such that for every prime number $p$ there
exists $\epsilon_p\in\wtN(\Qp)$ such that the conjugacy class of
$\epsilon_p\trho_p(\Fr_v)$ is $\Cl(\wtFr_v)$. 
\end{prop}
The proof requires several lemmas and the following notation.
The group $\Aut'$ defined above also acts on $G_A^\der$. 
As $\Out'(\wtG^{s,\der})\subset\Out(G_A^{s,\der})$ acts trivially on
the centre of $G_A^{s,\der}$, this action extends to an action of $\Aut'$
on $G_A$ with trivial action on the centre.
As above, we write $\Cl(G_A)$ for the
categorical quotient of $G_A$ by this action and let
$\Cl(\rho_{A,p})\colon\cG_F\rightarrow\Cl(G_A)(\Qp)$ 
be the map induced by $\rho_{A,p}\colon\cG_F\rightarrow G_A(\Qp)$, for
each prime number $p$. 

Similarly, $\Aut'$ acts naturally on $\wtG^\ad$ and on
$\wtG^\ab$ so we deduce an action on $G'=\wtG^\ad\times\wtG^\ab$.
Let $\Cl(\wtG^\ad)$ (resp.\ $\Cl(G')$) be the quotient of $\wtG^\ad$
(resp.\ $G'$) by $\Aut'$. 
The maps $G_A\rightarrow\wtG^\ad\rightarrow G'$ and 
$\wtG\rightarrow G'$ are $\Aut'$-equivariant and therefore induce maps 
$\Cl(G_A)\rightarrow\Cl(G')$ and $\Cl(\wtG)\rightarrow\Cl(G')$. 

Recall from~\ref{LADH_intro} and corollary~\ref{constr_LADH} that the
$\trho_p$ lift the representations
\[
\rho_p'=(\rho^\ad_p,\rho^\ab_p)\colon\cG_F\rightarrow G'(\Qp). 
\]
Let the $\Cl(\rho_p')\colon\cG_F\rightarrow\Cl(G')(\Qp)$ be the induced
maps. 
\begin{lemma}\label{classes_Fr-A}
For each $v\not\in\Sigma$, there is an element
\[
\Cl(\Fr_{A,v})\in\Cl(G_A)(\QQ)
\]
such that for every prime number $p$ with $v(p)=0$
we have 
$\Cl(\rho_p)(\Fr_v)=\Cl(\Fr_{A,v})$. 
\end{lemma}
\begin{proof}
We use the notation of~\ref{intro_char_pol}. 
The hypotheses that $G_A$ arises from the construction
of~\ref{better_unliftable_lift} imply that $G_A^\der\times L^\times$
acts on $\rH_\B^1(A(\CC),\QQ)$ and that this representation is a Weil
restriction of $W\otimes_{K_0}V^s$, where $W$ is the representation of
$L^\times$ on $L$ by left multiplication and $V^s$ is a multiple of
the representation of $G^{s,\der}$ of highest weight $\varpi_1$. 
It follows in particular that
$L$ lies in the centre of 
$\End^0(A_{/\CC})=\End(A_{/\CC})\otimes_\ZZ\QQ$.  

Since $v\not\in\Sigma$, the abelian variety $A$ has good
reduction $A_v$ at $v$ so we can identify $\End^0(A_{/\CC})$ with a
subalgebra of $\End^0(A_{\bar v})$.
Here $A_{\bar v}$ is the base extension of $A_v$ to the algebraic
closure of the residue field of $F$ at $v$. 
We obtain an embedding $L\subset\End^0(A_{\bar v})$. 
Let $\pi_v\colon A_v\rightarrow A_v$ be the Frobenius endomorphism.
For each prime number $p$ different from the residue characteristic
at $v$, there is a canonical, hence $L$-equivariant, identification 
$\rH_\et^1(A_{\bar F},\Qp)=\rH_\et(A_{\bar v},\Qp)$. 
Under this identification, the action of $\rho_p(\Fr_v)$ on the left hand
side corresponds to the action of $\pi_v$ on the right hand side. 

Since $\pi_v$ is semi-simple and lies in the centre of $\End^0(A_v)$,
the subalgebra $M=L[\pi_v]\subset\End^0(A_v)$ is a product of
number fields $M_i$. 
Each $M_i$ is of the form $\QQ(\alpha_i)$ for some
$\alpha_i\in\End^0(A_v)$ and it is a standard fact 
(cf.~\cite[\S19, theorem~4]{mumford:av}) that the characteristic
polynomial of $\alpha_i$ acting on $\rH_\et(A_{\bar v},\Qp)$ has
coefficients in $\QQ$ and is independent of $p$. 
This implies that there exists an $M$-module $U_v$ such that for every
prime number $p$ with $v(p)=0$ there is an isomorphism 
$\rH_\et(A_{\bar v},\Qp)\cong U_v\otimes_\QQ\Qp$ 
of $M\otimes_\QQ\Qp$-modules. 
Let $Q[X]\in L[X]$ be the characteristic polynomial of $\pi_v\in M$ acting
on $U_v$ as an $L$-linear endomorphism. 
For every $p$ with $v(p)=0$, the characteristic polynomial of
$\rho_{A,p}(\Fr_v)$, acting $L\otimes_\QQ\Qp$-linearly on 
$\rH_\et^1(A_{\bar F},\Qp)$, is equal to $Q$. 
This result is due to Shimura, see
\cite[11.10.1]{shimura:algebraic_n-fields}. 

Consider the map $G^{s,\der}_{/\bQ}\rightarrow\AA^n_{/\bQ}$
corresponding to the characteristic polynomial in the representation
with highest weight $\varpi_1$.
It factors through the quotient
of $G^{s,\der}$ for the action of $\Aut'(G^{s,\der})$, giving a
map $\Cl(G_A)\rightarrow\Res_{L/\QQ}\AA^n_{L}$. 
This map is easily seen to be injective on geometric points. 
The fact that $Q$ is the characteristic polynomial of $\rho_{A,p}(\Fr_v)$
for every $p$ with $v(p)=0$ implies that $Q$ is the image in
$\Res_{L/\QQ}\AA^n_{L}$ of an element $\Cl(\Fr_{A,v})\in\Cl(G_A)(\QQ)$ and
that this element $\Cl(\Fr_{A,v})$ verifies the condition of the lemma. 
\end{proof}
\begin{cor}\label{G'-class_rational}
For each $v\not\in\Sigma$, there is an element
\[
\Cl(\Fr'_v)\in\Cl(G')(\QQ)
\]
such that for every prime number $p$ with $v(p)=0$
we have 
$\Cl(\rho'_p)(\Fr_v)=\Cl(\Fr'_v)$. 
\end{cor}
\begin{proof}
By \ref{LADH_intro}, the
$\rho^\ab_p\colon\cG_F\rightarrow\wtG^\ab(\Qp)$ form a compatible
system, so all $\rho^\ab_p(\Fr_v)$ are defined over $\QQ$ and
coincide. 
This gives rise to an element $\Fr_v^\ab\in\wtG^\ab(\QQ)$. 
On the adjoint side, let $\Cl(\wtFr^\ad_v)\in\Cl(\wtG^\ad)(\QQ)$ be
the image of $\Cl(\Fr_{A,v})$ and let $\Cl(\Fr'_v)$ be the image of 
$(\Cl(\wtFr^\ad_v),\Fr_v^\ab)$ in $\Cl(G')(\QQ)$. 

For each $p$, the representation $\trho_p$ is a weak geometric lift of
$\rho_{A,p}$ so the projection 
$\trho_p^\ad\colon\cG_F\rightarrow\wtG^\ad(\Qp)$ coincides with the
composite of $\rho_{A,p}$ with the projection
$G_A(\Qp)\rightarrow\wtG^\ad(\Qp)$. 
It follows that $\Cl(\Fr'_v)$ is the element promised by
the corollary. 
\end{proof}
\begin{lemma}\label{Hasse_princ}
Let $M$ be either the group scheme $\wtN=\Res_{K_0/\QQ}\wtN^s$
(the centre of $\wtG^\der$, cf.~\ref{intro_char_pol}) or 
\[
\Out'(\wtG)=\Aut'(\wtG)/\wtG^\ad\cong\Res_{K_0/\QQ}\mu_2
\]
and let $P$ be a co-finite set of prime numbers. 
Then the natural map 
\[
\rH^1_\et(\Spec(\QQ),M)\rightarrow
\prod_{p\in P}\rH^1_\et(\Spec(\Qp),M)
\]
is injective.
\end{lemma}
\begin{proof}
In the lemma, $M$ denotes either $\wtN$ or $\Out'(\wtG)$. 
We put $M'=\wtN^s$ in the former and $M'=\mu_2$ in the latter case.
Let $f\colon\Spec(K_0)\rightarrow\Spec(\QQ)$ be the natural morphism. 
This implies that $M=f_*M'$ as \'etale sheaves and since $f_*$ is exact by 
\cite[{I}{I},~Corollary~3.6]{milne:etale}, this gives an isomorphism
\[
\rH^1_\et(\Spec(\QQ),M)\cong\rH^1_\et(\Spec(K_0),M'),
\]
cf.~\cite[I.~2.5]{serre:cohom_gal5} for the interpretation in terms of
Galois cohomology. 
Similarly, for every prime number $p$, there is an isomorphism 
\[
\rH^1_\et(\Spec(\Qp),M)\cong
\prod_{\gp\mid p}\rH^1_\et(\Spec(K_{0,\gp}),M'),
\]
obtained from the above identification by base change to $\Spec(\Qp)$. 
To prove the lemma, it is sufficient to prove the injectivity of the
map 
\[
\rH^1_\et(\Spec(K_0),M')\rightarrow
\prod_{\gp\in P'}\rH^1_\et(\Spec(K_{0,\gp}),M'),
\]
where $P'$ is the set of primes of $K_0$ lying over the rational
primes in $P$.
In the case where $M=\Out'(\wtG)$, the group scheme $M'=\mu_2$ is
constant and as $\rH_\et^1(\Spec(K),\mu_2)\cong\Hom(\cG_K,\mu_2(K))$
for any field $K$, the lemma follows. 

In the case where $M'\cong\wtN^s$, there are two possibilities, 
$M'$ is geometrically isomorphic to either $(\ZZ/2\ZZ)^2$ or
$\ZZ/4\ZZ$. 
In both cases the group scheme becomes trivial over an
extension of degree (at most) $2$. 
If $M'$ is already trivial over $K_0$, then the above argument
applies. 
If $M'$ is geometrically isomorphic to $(\ZZ/2\ZZ)^2$ and trivial over
the quadratic extension $K\supset K_0$, then $M'$ is a Weil
restriction and the lemma follows by the same argument as before. 

We are left with the case where $M'_K\cong(\ZZ/4\ZZ)_K$ for a quadratic
extension $K\supset K_0$. 
The proof in this remaining case, which also applies in the other cases where
$M'\cong\wtN^s$, makes use of the long exact cohomology sequences
associated to the short exact sequence 
$1\rightarrow\mu_2\rightarrow M'\rightarrow\mu_2\rightarrow1$ of
\'etale sheaves on $\Spec(K_0)$. 
This gives a commutative diagram with exact rows,  
\[
\begin{CD}
\mu_2(K_0)@>>>\rH^1(K_0,\mu_2)@>i>>
\rH^1(K_0,M')@>>>\rH^1(K_0,\mu_2)\\
@VVV @VV\loc_{\mu_2}V @VV\loc_{M'}V @VV\loc_{\mu_2}V\\
\prod\mu_2(K_{0,\gp})@>>>\prod\rH^1(K_{0,\gp},\mu_2)@>\prod i_\gp>>
\prod\rH^1(K_{0,\gp},M')@>>>\prod\rH^1(K_{0,\gp},\mu_2).
\end{CD}
\]
Here all $\rH^1$ are \'etale cohomology groups over the spectrum of
the specified field. 
These groups can be
identified with the corresponding Galois cohomology groups. 
The products in the second row are over all $\gp\in P'$. 

The image of $i$ identifies with 
$K_0^\times/\langle K_0^{\times2},\alpha\rangle$ for some 
$\alpha\in K_0^\times$, where $\langle K_0^{\times2},\alpha\rangle$ is
the subgroup of $K_0^\times$ generated by $\alpha$ and the squares in
$K_0^\times$. 
Similarly, the image of each $i_\gp$ is identified with 
$K_{0,\gp}^\times/\langle K_{0,\gp}^{\times2},\alpha\rangle$.
If $x\in\ker\loc_{M'}$, then the injectivity of $\loc_{\mu_2}$ implies
that $x$ lies in the image of $i$. 
Let $y\in K_0^\times$ represent the preimage of $x$ in
$K_0^\times/\langle K_0^{\times2},\alpha\rangle$. 
Since $x\in\ker\loc_{M'}$, the element $y$ lies in 
$\langle K_{0,\gp}^{\times2},\alpha\rangle$ for each $\gp\in P'$, 
so $y$ is a square in $K_0(\sqrt\alpha)$ locally at each place
above a $\gp\in P'$. 
It follows that $y$ is a square in $K_0(\sqrt\alpha)$, hence it maps to
$1$ in $K_0^\times/\langle K_0^{\times2},\alpha\rangle$.
\end{proof}
\begin{proof}[Proof of proposition~\ref{classes_Fr}]
The map $\wtG\rightarrow G'$ induces a map
$\pr\colon\Cl(\wtG)\rightarrow\Cl(G')$. 
Obviously,
\[
\Cl(\rho'_p)=\pr_{\Qp}\circ\Cl(\trho_p)\colon\cG_F\rightarrow\Cl(G')(\Qp)
\]
for each $p$. 

For the groups $\wtG$ and $G'$, let $\Conj(\wtG)$ and $\Conj(G')$ be
the varieties of geometric conjugacy classes, i.\ e.\ the quotients
for the actions of $\wtG^\ad$. 
Note that $\Cl(\wtG)$ (resp.\ $\Cl(G')$) is the quotient of
$\Conj(\wtG)$ (resp.\ $\Conj(G')$) by $\Out'(\wtG)$ and that 
there is a commutative diagram 
\[
\begin{CD}
\Conj(\wtG)@>>>\Cl(\wtG)\\
@VVV @VV\pr V\\
\Conj(G')@>>>\Cl(G').
\end{CD}
\]
There is a Zariski closed subset $B\subset\Cl(G')$ such that 
$\Conj(G')\rightarrow\Cl(G')$ and 
$\pr$ are both unramified outside $B$. 
Let the $\Cl(\Fr'_v)$ be as in corollary~\ref{G'-class_rational}.
The hypothesis on the rank of the Zariski closure of the image of
$\rho_{A,\ell}$ and the Chebotatev density theorem imply that 
there is a
set $\Sigma_{\Fr}$ of places of $F$ of Dirichlet density $0$ such that 
$\Cl(\Fr'_v)\in B(\QQ)$ if and only if $v\in\Sigma_{\Fr}$. 
To prove the proposition, it suffices to show that for every 
$v\not\in\Sigma_{\Fr}$, the
fibre of $\pr$ over $\Cl(\Fr'_v)$ contains an element 
$\Cl(\wtFr_v)\in\Cl(\wtG)(\QQ)$. 

Fix $v\not\in\Sigma_{\Fr}$. 
The fibre of $\Conj(G')\rightarrow\Cl(G')$ over $\Cl(\Fr'_v)$ is a
$\Out'$-torsor which admits a $\Qp$-valued point
$\Conj(\rho_p'(\Fr_v))$ for almost every $p$. 
It follows from lemma~\ref{Hasse_princ} that there exists
$\Conj(\Fr'_v)\in\Conj(G')(\QQ)$ above $\Cl(\Fr'_v)$. 
The fibre of $\Conj(\wtG)\rightarrow\Conj(G')$ over this element is a
$\wtN$-torsor. 
As this fibre identifies with $\pr^{-1}(\Cl(\Fr_v'))$, this
last fibre is also a $\wtN$-torsor. 
It has a $\Qp$-valued point $\Cl(\trho_p(\Fr_v))$ for almost
every $p$, so the promised existence of
$\Cl(\wtFr_v)\in\Cl(\wtG)(\QQ)$ follows 
from another application of lemma~\ref{Hasse_princ}.
\end{proof}
\addcontentsline{toc}{section}{References}
\bibliography{biblio}
\bibliographystyle{mybib}
\setbox1=\hbox{\tt http://www-irma.u-strasbg.fr/\~{}noot}
\noindent
Rutger Noot\\
I.~R.~M.~A.\\
Universit\'e Louis Pasteur and CNRS\\
7 rue Ren\'e Descartes\\
67084 Strasbourg, France\\
{\tt rutger.noot@math.u-strasbg.fr}\\ 
{\tt http://www-irma.u-strasbg.fr/\~{}noot}
\end{document}